\documentstyle[amscd,amssymb,bbm,verbatim]{article}
 
\newtheorem{thm}{Theorem}[section]

\newtheorem{cor}{Corollary}[section]
\newtheorem{prop}{Proposition}[section]

\newtheorem{lemma}{Lemma}[section]

\newcounter{num}
\begin{document}

\font\bbb = msbm10
\def\Bbb#1{\hbox{\bbb #1}}
\font\lbbb = msbm7
\def\LBbb#1{\hbox{\lbbb #1}}
\def\squarebox#1{\hbox to #1{\hfill\vbox to #1{\vfill}}}
\newcommand{\qed}{\hspace*{\fill}
           \vbox{\hrule\hbox{\vrule\squarebox{.667em}\vrule}\hrule}\smallskip}

\newcommand{\beq}        {\begin{eqnarray}}
\newcommand{\eeq}        {\end{eqnarray}}
\newcommand{\be}         {\begin{equation}}  
\newcommand{\ee}         {\end{equation}}
\newcommand{\beqn}       {\begin{eqnarray*}}  
\newcommand{\eeqn}       {\end{eqnarray*}}

\def\ep{\varepsilon}
\def\R{\mathbb R}
\def\C{\mathbb C}
\def\a{{\alpha}}
\def\b{{\beta}}
\def\p{\partial}
\def\t{{\theta}}
\def\T{\mathbb T}
\def\C{{\mathbb C}}
\def\Z{{\mathbb Z}}
\def\I{{\mathcal I}}
\def\E{{\mathbb E}}
\def\l{{\ell}}
\def\ld{{\lambda}}
\def\G{{\Gamma}}
\def\d{{\delta}}
\def\pr{{\prime}}
\def\up{{\uparrow}}
\def\ra{{\rightarrow}}
\def\m{{\mu}}
\def\mm{{\mathfrak{m}}}
\def\sn{{\mbox{sn}}}
\def\cn{{\mbox{cn}}}
\def\dn{{\mbox{dn}}}
\def\sech{{\mbox{sech}}}
\def\g{{\mathfrak{g}}}

\title{On the shape of the ground state eigenvalue density
of a random Hill's equation}

\date{November 11, 2004}

\author{ Santiago Cambronero${\ }^{*}$
\and Jos{\'e} Ram{\'\i}rez 
\thanks{ Department of Mathematics, Universidad de Costa Rica, San Jose
2060, Costa Rica.
scambro@emate.ucr.ac.cr, jramirez@emate.ucr.ac.cr}
\and Brian Rider
\thanks{Department of
Mathematics,
University of Colorado at Boulder, 
UCB 395, Boulder CO 80309. brider@euclid.colorado.edu} 
}

\maketitle

\begin{abstract} 
Consider the Hill's operator $Q = - d^2/dx^2 + q(x)$ in which 
$q(x)$, $0 \le x \le 1$, is a White Noise.  Denote
by $f(\m)$ the probability density function of $-\ld_0(q)$,
the negative of the ground state eigenvalue, at $\m$.  We
prove the detailed asymptotics:
\[
   f(\m) 
= \frac{4}{3 \pi} \, \m \, \exp \Bigl[ - \frac{8}{3} \m^{3/2} 
- \frac{1}{2} \m^{1/2}
   \Bigr] ( 1 + o(1) )      
\]
as $\m \ra +\infty$.  This result is based on a precise Laplace
analysis of a functional integral representation for $f(\m)$
established by S. Cambronero and H.P. McKean in \cite{CM}.
\end{abstract}

\section{Introduction}
\setcounter{num}{1}
\setcounter{equation}{0}
\label{sec:int}

We consider fluctuations of the ground state eigenvalue 
of a random Hill's operator $Q = -d^2/dx^2 + q(x)$
in which the periodicity is fixed at one and the potential $q(x)$ 
is a White Noise. Formally, $q(x) = b^{\prime}(x)$
for a standard Brownian Motion $b(x)$, and the
eigenvalue problem $Q \psi = \ld \psi$  must be interpreted in the
sense of It$\hat{\mbox{o}}$. In those terms it reads 
$d \psi^{\prime}(x) = d b(x) - \ld \psi(x) dx$,  and
there is no problem solving for $\psi$ in the space $C^{3/2-}$.  

Let $\ld_0(q)$ denote the ground state eigenvalue of $Q$.
Our jumping off point is an explicit formula for the law of
this object due to the authors of \cite{CM}.  
Choosing the sign for later convenience we further denote by $f(\m)$
the probability density function of $- \ld_0(q)$  at the point $\m$.
Then, the result of \cite{CM} is
\be
\label{equ:dens}
          f(\m)
          =   \frac{1}{\sqrt{2 \pi}} \int_{H} \exp{ \Bigl\{ -\frac{1}{2} 
              \int_0^1 (\mu - p^2(x))^2 dx  \Bigr\}} A(p) \, d P_0
\ee
where 
\be
\label{equ:A}
   A(p) =  \int_0^1 e^{2 \int_0^x p(x^{\pr}) d x^{\pr}} dx  \times 
                 \int_0^1 e^{- 2 \int_0^x p(x^{\pr}) d x^{\pr}} dx,
\ee
and $P_0$ is a probability measure on
$H$, the space of continuous functions of period one and mean zero.
More specifically, ${P}_0$ is the  
Circular Brownian Motion ($CBM$) $p(x)$, $ 0 \le x \le 1$,
conditioned so that $\int_0^1 p(x) dx = 0$.   $CBM$, we recall, is the measure 
on periodic paths formed from the standard Brownian Motion starting from 
$p(0) = c$, conditioned to return to $c$ at 
$x = 1$, with this common starting/ending point distributed
over the line according to Lebesgue measure. In other words, for
any event ${\cal A}$ of the path, $CBM({\cal A}) = \frac{1}{\sqrt{2 \pi}} 
\int_{-\infty}^{\infty}
BM_{00} ({\cal A} + c) dc$ in which $BM_{00}$ denotes the mean of the Brownian 
Bridge of length one.  While $CBM$ itself has infinite total mass,
$P_0$ is the distribution of an honest rotation invariant Gaussian process: 
a short computation will show that
\[
  {E}_0 \Bigl[ F(p) \Bigr] = 
                     BM_{00} \Bigl[  F ( p - \int_0^1 {p} ) \Bigr]
\]
for any bounded measurable test function $F = F(p(x): 0 \le x \le 1)$.

The functional integral representation of the density 
$f(\m)$ given by (\ref{equ:dens}) is based on a correspondence
between Hill's equation and Ricatti's equation.  To explain, first    
fix a $\ld$ to the left of $\ld_0(q)$.  Then it is well known (see \cite{MW})
that Hill's equation possesses a positive solution $\psi$ with multiplier $m > 0$: 
$-\psi^{\prime \prime} + q \psi =  \ld \psi$ and $\psi(x+1) = m \psi(x)$.
Further, the logarithmic derivative $p = \psi^{\prime}/\psi$ is a periodic solution 
of Ricatti's equation $ q = \ld + p^{\prime} + p^2$ with $\int_0^1 p \ge 0$. 
Now view this {\em Ricatti correspondence} as a map between measure spaces. 
Under this map the restriction of the White Noise 
measure to $\{  \ld \le \ld_0(q) \}$ is identified 
with the restriction of the $CBM$ to $\{ \int_0^1 p \ge 0 \} $ up to a suitable 
Jacobian factor.  The computation of that Jacobian is the chief accomplishment of 
\cite{CM},  resulting in, among other things, the formula (\ref{equ:dens}).

Given (\ref{equ:dens}), 
our purpose here is to describe the shape of the density.
That is, we investigate the detailed asymptotics of $f(\m)$ as $\m \ra \pm \infty$.
The left tail has in fact already been discussed in
\cite{CM}: they remark that
\[
          f(\mu) = \sqrt{\frac{ (-\m) }{\pi}} \exp{ \Bigl[  - \frac{1}{2} \m^2 
          - \frac{1}{\sqrt{2}} (-\m)^{1/2} \Bigr]} (1 +o(1))
           \;\;\; \mbox{ for } \m \ra -\infty.
\]
This is easy to understand.  The exponent is expanded out as in
$ \int_0^1 (|\m| + p^2)^2 =  \m^2 - 2 |\m| \int_0^1 p^2 - \int_0^1 p^4$.
Next, $A(p)$ and $\exp \{ -(1/2) \int_0^1 p^4 \}$ are shown to be negligible 
when compared to
the coercive Gaussian weight $ \exp \{ - |\m| \int_0^1 p^2 \}$. 
The derivation is completed by 
computing 
$ E_0 [  \exp \{ -|\m| \int_0^1 p^2 \} ] $
exactly.

The analysis of right tail is far more involved. Our result is
the following.

\begin{thm}
The probability density of $-\ld_0(q)$ has the
shape
\label{thm:main}
\be
\label{equ:main}
  f(\mu) =   \frac{4}{3 \pi}  \, \m \,  
  \exp{\Bigl[ - \frac{8}{3} \m^{3/2} - \frac{1}{2} \m^{1/2} \Bigr]}  ( 1 + o(1) )
\; \; \;  {\rm{ for }}  \; \; \m \ra + \infty.
\ee
\end{thm}

While the limit $\m \ra -\infty$ concentrates the path at the unique and trivial trajectory
$p \equiv 0$, taking $\m \ra +\infty$  leads to a completely different picture.  In this
regime it is most advantageous for $p$ 
to live near $\pm \sqrt{\m}$.  However, because of the  restriction to $\int_0^1 p = 0$, 
the path is forced to divide its favors more or less evenly between the two choices 
of sign.  Furthermore, due to the rotation 
invariance of CBM, any  translation of an extremal path is also extremal.
Thus, for $\m \ra +\infty$, we are dealing with Laplace asymptotics of 
a degenerate function space integral.  Even so, we are able to obtain beyond
leading order information.

There is of course no shortage of investigations into the precise large deviations 
of Wiener-type or other functional integrals, including cases in with the underlying 
functional possesses degeneracies.  Important examples are \cite{ER}, \cite{Bol},
\cite{GB}, \cite{KS1} and \cite{KS2}.  Nevertheless, the present problem has various
features which set it apart and require the analysis to be done ``by hand''.
First, the large parameter enters $f(\m)$ in a fundamentally different way
than is assumed throughout the cited list.
More importantly, those dealing with degenerate problems 
assume nondegeneracy in directions orthogonal to the extremal set.  Our integral
does not posses this property; there 
exists a more subtle degeneracy besides that stemming from the translation invariance.

Random Schr\"odinger operators of the type 
$Q$ arise in models of disordered  solids as is 
explained in the comprehensive book \cite{LGP}.  The White Noise potential
offers a simplifying caricature.  Its use goes back to \cite{FL}, but
see also  \cite{Hlp} and \cite{FM} 
which discuss the integrated density of states.  A description 
of the ground state energy is of separate importance.  In the present 
White Noise setting with dimension equal to one, \cite{M1}
proves a limit law for $-\ld_0(q)$ as the periodicity is taken to infinity.
Also in the thermodynamic regime, the study of the almost sure behavior of the ground
state in any dimension subject to a potential of Poisson or Gibbs type is well 
developed: see \cite{Sn}  and \cite{Mrk} and references therein.  Still,
the understanding of the actual distribution 
in a finite volume with any kind of potential 
remains in its infancy, and the result described above
prompts further inquiry.  In particular, it is reasonable to ask to what
extent the shape of $f(\m)$ is universal for some class of rough potentials.
The Gaussian/sub-Gaussian tails to the left/right 
seen here have an intuitive explanation: level repulsion holding down the 
left tail, with the ground state free to take advantage of deep wells
created by the White Noise to the right. 

The rest of the paper takes the following course.
In Section 2 we study the associated rate function 
and  discuss a leading order result of the form 
$ \m^{-3/2} \log f(\m) \simeq -8/3$ for $\m \uparrow \infty$.  Asymptotics 
at this level are accounted for by a vicinity of a one parameter family of 
paths (the degeneracy).   
Section 3 outlines how to expand about the set of extrema.
Namely, the degeneracy is dealt with by a conditioning procedure, leaving 
an integral with respect to a certain Gaussian measure as the principle lower
order term.  Using a connection to a particular (deterministic) Hill operator, 
required properties of this Gaussian measure are collected in Section 4. 
The calculation is picked up again in
Section 5 which contains the main error estimate: here we dispose of the terms beyond the
Gaussian correction.
Afterward, in Section 6, the Gaussian correction is computed exactly.
In essence this completes our calculation.  Section 7 gathers the results through
that point and states Theorem \ref{thm:mainish} which, in some sense, is a more accurate statement of
the main result.   It is then explained how asymptotics of the various quantities appearing
in Theorem \ref{thm:mainish} translate that result into the above Theorem \ref{thm:main}.  Finally,
Section 8 serves as an appendix containing various technicalities  needed along the way.

\section{Leading order asymptotics}
\setcounter{num}{2}
\setcounter{equation}{0}
\label{sec:LD}

\subsection{The rate function}
\label{sec:rate}

At an exponential scale, the asymptotics of integrals of the type
(\ref{equ:dens}) are well understood to be associated with a 
characteristic  variational problem. 
In the case of $f(\mu)$, that problem is to minimize
\be
\label{equ:Imu}
   I_{\mu}(p) \equiv \frac{1}{2} \int_0^1 \Bigl(\mu - p^2(x) \Bigr)^2 dx  
           + \frac{1}{2} \int_0^1 ( p^{\pr}(x) )^2 dx
\ee
over periodic functions of mean zero.  Scaling
as in $p(\cdot) \ra \frac{1}{\sqrt{\m}}
f( \cdot /\sqrt{\m} )$ we find that
\beq
\label{equ:scaledvar}
  {\inf}_{p \in H} I_{\mu} (p) 
  & = & 
  \m^{3/2} \ {\inf}_{f \in H_a} \Bigl\{ \frac{1}{2} \int_{-a}^a (1 - f^2(x))^2 dx + \frac{1}{2} 
                           \int_{-a}^{a} (f^{\pr}(x))^2 dx \Bigr\}  \\
  & \equiv & \m^{3/2} \ {\inf}_{f \in H_a} I(f; a) \nonumber
\eeq
in which $ a = \sqrt{\m}/2$ and now
$H_a$ is the class of periodic $C^{1}$ functions satisfying
$\int_{-a}^{a} f = 0$.
This already suggests the $3/2$-power in the exponent of (\ref{equ:main}).  
As to its $8/3$-multiplier and further properties of the {\em rate function} $I$ we have
the following.
 
\begin{thm}  
\label{thm:rate} The infimum 
$I^*(a) = \inf_{f \in H_a} I(f; a)$ is attained.  Further, 
$I^*(a) \le 8/3$ for all $a > 0$
and $ I^*(a) \ra 8/3$ as $a \ra \infty$.   
\end{thm}

\noindent{\bf Proof}  As the level sets $\{ I(\cdot; a) \le K \}$ are weakly compact in $H^1$
and strongly so in $L^{\infty}$, the existence of a minimizer poses no problem.  For the upper
bound on $I^*(a)$, consider the test function
\be
\label{equ:test}
    f_a(x) =\left\{ 
\begin{array}{lll}
-\tanh {(x+a)} & \mbox{for} & -a\le x<-{a}/{2}, \\ 
\tanh {(x)} & \mbox{for} & -{a}/{2}\le x< {a}/{2}, \\ 
-\tanh {(x-a)} & \mbox{for} & {a}/{2}\le x<a.
\end{array}
\right. 
\ee
Then,
\beqn
\label{equ:test1}
  I^*(a) \le I(f_a; a)  & = &   \int_{-a/2}^{a/2} (1 - f_a^2(x) )^2 dx 
                          +   \int_{-a/2}^{a/2} (f_a^{\pr}(x))^2  \\
& \le &    \int_{-\infty}^{\infty} (1 - \tanh^2(x) )^2 dx 
                          +   \int_{-\infty}^{\infty} (\tanh^{\pr}(x))^2 dx  \nonumber \\
                   &  = & 2 \int_{-\infty}^{\infty} {\sech}^4(x) dx = {8}/{3}. \nonumber
\eeqn

As for the convergence $I^*(a) \ra 8/3$, computing a first variation of $I(\cdot; a)$
shows that any minimizer satisfies $f^{\pr \pr} = 2 f^3 - 2 f^2 - \a$ where
$\a$ is the Lagrange multiplier corresponding to the condition $\int_{-a}^{a} f = 0$. Multiplying
by $f^{\pr}$ and integrating, this may be brought into the form
\be
\label{equ:EL1}
   \frac{1}{2} (f^{\pr})^2 = \frac{1}{2} f^4 - f^2 - \a f + \frac{1}{2} \b
\ee
with a new constant $\b$.  We will show in the appendix that, as $a \ra \infty$, 
$\a \ra 0$ and $\b \ra 1$. 

Now choose a minimizer $f_a^*$ for each $a \gg 1$.  Since every $f_a^*$ must have at least
two roots on account of being mean zero, by rotation invariance we can assume that 
$f_a^*(0) = 0$ for all $\{f_a^*\}$. Next, for any large fixed $M > 0 $ less
than $a/2$ we certainly have that
\[
  I^*(a)   =  \frac{1}{2} \int_{-a}^{a} (d f_a^*/dx )^2  + \frac{1}{2} \int_{-a}^{a} (1 - (f_a^*)^2)^2 
           \ge  \frac{1}{2}\int_{-M}^{M} (d f_a^*/ dx)^2  + \frac{1}{2} \int_{-M}^{M} (1 - (f_a^*)^2)^2. 
\]
It follows that the sequence $\{ f_a^*\}$
is bounded in $H^1 \cap L^{\infty} [ -M, M]$ and so has a subsequence converging weakly in $H^1$
and strongly in $L^4$.  By writing the relation (\ref{equ:EL1}) in the form
\be
f^{\pr} =  \sqrt{f^4 - 2f^2 - 2\a f + \b},
\ee
using its weak formulation and the regularity theory that comes with it, we
find that the above 
convergence is sufficient to conclude that any
limit satisfies $f_{\infty}^{\pr} = 1 - f_{\infty}^2$ over $[-M, M]$ 
and $f_{\infty}(0) = 0$.  That is, 
$f_{\infty} = \tanh$. 

Since the centered $\{ f_a^{*} \}$ converge to $f_{\infty}$ 
on an interval $[-M, M]$ for any
choice of $M > 0$, it follows that  the distance between any two zeros of $f_a^*$
must be tending to infinity as $a \ra \infty$.  Therefore, we may also
isolate a symmetric interval of length $2M$ about a second zero of $f_a^*$.
On this second interval we will
have the same type of convergence (by precisely the same arguments). Then,
by adding both contributions we also conclude that 
\[
  I^*(a)    \ge \int_{-M}^{M} (d f_a^*/ dx)^2  + \int_{-M}^{M} (1 - (f_a^*)^2)^2
\]
for all $a$ large enough.
This results in 
\[
\liminf_{a \ra \infty} I^*(a)
\ge 2 \int_{-M}^M \sech^4(x) dx  \uparrow 8/3
\]
upon letting $M \up \infty$ afterward. 
The proof is finished.

\subsection{Large Deviations}
\label{sec:rare}

While the previous result serves as a guide, the next step is to extract the $\exp{[-8/3 \mu^{3/2}]}$
behavior from the integral $f(\m)$.

\begin{thm} 
\label{thm:LD}
We have the leading order asymptotics: 
\[
   \limsup_{\m \ra \infty} \m^{- 3/2}  \log  \int_{H } 
   e^{-\frac{1}{2} \int_0^1 (\m - p^2)^2 } A(p) \ d P_0  \le  - {8}/{3}.
\]
\end{thm}

More important in the sequel, we show that one has sharper decay of the same order
when the integral is restricted to a set away from any $I_{\m}$ minimizer.  
Define
\be
\label{equ:nearmin}
           C_{\ep}^{\m} = \Bigl\{ p \in H : d(p, {\cal M}) \le \ep \sqrt{\m} \Bigr\}
\ee
in which ${\cal M}$ is the set of minimizers of $I_{\mu}$ and 
$d(p, {\cal A})$ is the distance between a path $p$ and a set ${\cal A}$ 
in sup-norm.  Then we have:
 
\begin{thm} 
\label{thm:rare}
There exists a $\eta > 0$ depending on $\ep$ 
and a constant $C$ so that
\[
   \int_{ H \backslash C_{\ep}^{\mu}}  
                e^{-\frac{1}{2} \int_0^1 (\m - p^2)^2 } A(p) \ d P_0
   \le {C} e^{- (8/3 + \eta) \mu^{3/2}}
\]
for all $\m$ large enough.
\end{thm}

\bigskip

\noindent{\bf Proof of Theorem {\ref{thm:LD}}} 
Matters are simplified by noticing 
that it is enough to prove that
\be
\label{sinA}
 \limsup_{\mu \ra \infty } \mu^{-3/2} \log \int_{H} e^{-
 \frac{1}{2}\int_{0}^{1}\left( \mu - p^{2}\right) ^{2}} 
 dP_{0}  \le - \frac{8}{3}.
\ee
This is because 
$A( p ) \leq \exp[ 2\sqrt{\int_{0}^{1}p^{2}} ]
\leq \exp  [ 1+\int_{0}^{1}p^{2} ]$
which implies the upper bound
\[
\int_H e^{-\frac{1}{2}\int_{0}^{1}\left( \mu -p^{2}\right) ^{2}}A\left(
p\right) d P_{0}\leq e^{\mu +3/2}\int_H  e^{-\frac{1}{2}\int_{0}^{1}\left( \mu
+1-p^{2}\right)^{2}} d  P_{0},
\]
the prefactor $e^{\mu + 3/2}$ being irrelevant in the present scale.

The proof of (\ref{sinA}) is split into several steps.  
First we define the set of paths
\be
  H^{\mu } ( \gamma ,\eta ) = \left \{ p \in H:\left| \frac{1}{
2}\int_{0}^{1}\left( \mu -p^{2}\right)^{2}-\gamma \m^{3/2}\right| <\eta
\m^{3/2}\right\}
\label{DefH}
\ee 
for any positive $\eta$ and $\gamma$.  As we shall see, restricted to
such a set, the $P_0$ integral of $\exp{[-(1/2) \int_0^1 (\mu - p^2)^2]}$ 
is easy to control through the variational problem studied Section
{\ref{sec:rate}}.  
For this reason we make the decomposition 
\be
\label{equ:decomp}
\int_H e^{-\frac{1}{2}\int_{0}^{1}\left( \mu -p^{2}\right)^{2}} dP_{0} 
\le
\sum_{0 \le k \le {8}/{3}\eta } \int_{H^{\mu}(\eta k, \eta)}
e^{-\frac{1}{2} \int_{0}^{1} \left( \mu - p^{2}\right)^{2}} dP_{0} + \exp
\left( - \frac{8}{3} \mu^{3/2}\right)  
\ee 
after which we may invoke the bound 
\be
\label{equ:hd1}
\int_{H^{\m}(\gamma, \eta)} e^{-\frac{1}{2}\int_{0}^{1}\left( \mu
-p^{2}\right) ^{2}} d {P}_{0}\leq \exp \left[ \left( -\gamma +\eta \right)
\mu ^{3/2}\right] P_{0} \Bigl( H^{\mu }(\gamma, \eta) \Bigr)
\ee 
in each integral on the right hand side of (\ref{equ:decomp}).
This follows directly from the definition of $H^{\m}(\gamma, \eta)$.

The upshot is that we must now control  the probabilities $P_{0} ( H^{\mu
}(\gamma, \eta) )$.  Toward this end we overestimate further as in 
\beq
\label{equ:hd}
P_{0}\Bigl( H^{\m}(\gamma, \eta) \Bigr) & \leq & P_{0} \Bigl( \frac{1}{2}
\int_{0}^{1}\left( \mu -p^{2}\right)^{2}\leq \left( \gamma +\eta \right)
\mu ^{3/2}\Bigr) \\  & \equiv & P_{0}^{\m} \Bigl(
\frac{1}{2}\int_{0}^{1}\left( 1-p^{2}\right) ^{2}\leq \left( \gamma +\eta
\right) \frac{1}{\sqrt{\m}} \Bigr) \nonumber  
\eeq 
where we have
introduced the scaled measure  $P_{0}^{\m }( p \in A) \equiv P_{0} (
\frac{1}{\sqrt{\mu}} p \in A )$.  Under this scaling,  the previous display
reads   $P_{0}( H^{\mu}(\gamma, \eta ) \, ) \le P_{0}^{\mu }(
D^{\mu}(\gamma,\eta )\,)$, and  
\be
\label{equ:defD}
   D^{\mu} =   D^{\mu}(\gamma ,\eta) \equiv \left\{ p \in H:\frac{1}{2}
\int_{0}^{1}(1-p^{2})^{2}\leq \left( \gamma +\eta \right)
\frac{1}{\sqrt{\m}}  \right\}  
\ee 
marks yet another definition.

The next step, estimating the $P_0^{\mu}$ probability of $D^{\mu}$, is
accomplished by discretizing the path.  This is a common procedure, see for
example the proof of Schilder's Theorem in \cite{Var}.  Let $p_{n}$ be the polygonal path determined by the
values $p(k/n)$ at  $k/n$ for $k=0,\dots, n$ and introduce
$\widehat{p}_{n}= p_{n}- \int_{0}^{1} p_{n}$ to force things to reside in
$H$.   Now, for whatever set $D \subseteq H$ we have that, if $p\in D$,
then the polygonal  $p_n$ is either very close to $p$ or far away from $D$.
In symbols 
\be
\label{equ:enlarge}
  P_0^{\mu} (D) \le P_0^{\mu} \Bigl( || p - \widehat{p}_{n} || _{\infty }
\ge \delta \Bigr) + P_0^{\mu}  \Bigl(  \widehat{p}_n \in D_{\delta} \Bigr)
\ee 
in which $||p|| _{\infty } = \sup_{0 \le x < 1} |p(x)|$ for any  path $ p
\in H$ and $D_{\delta}$ is the $\delta$-enlargement of $D$ in that norm.
That is, $q \in D_{\delta}$ when  $ \inf_{p \in D} || q - p ||_{\infty} \le
\delta$.   We next tackle each term of the right of (\ref{equ:enlarge}).

For the deviation between $p$ and its approximate $\widehat{p}_n$ we first
note
\[
P_{0}^{\m} \Bigl( \left\| p-\widehat{p}_{n} \right\| _{\infty } \ge \delta
\Bigr) \leq P_{0}^{\m} \Bigl( \left\| p-p_{n}\right\| _{\infty }\ge
\frac{\delta }{2} \Bigr) ,
\]
because  $ \left\| p - \widehat{p}_{n}\right\| _{\infty }\leq \left\|
p-p_{n}\right\| _{\infty }+ | \int_{0}^{1} \left( p- p_{n}\right)  | \leq
2\left\| p - p_{n}\right\| _{\infty }. $ Then, 
\beq
\label{equ:en0}
 P_{0}^{\m} \Bigl( \left\| p- p_{n}\right\| _{\infty } \ge {\delta}/{2}
\Bigr)     & \le &  P_{0}^{\m} \Bigl( {\bigcup}_{k = 0}^{n-1} \Bigl\{
\sup_{ \frac{k}{n} \leq x \leq \frac{k+1}{n}}   \left| p(x)-
p({k}/{n})\right| \ge  {\delta }/{4} \Bigr\} \Bigr)  \\  & \le &  n
P_{0}^{\m} \Bigl( \sup_{0 \leq  x \leq 1/n} \left| p(x)- p(0) \right| \ge
{\delta }/{4} \Bigr),   \nonumber   
\eeq 
having used 
the rotation invariance of
$CBM$ in line two. Next we recall the definition of the measure
$P_{0}^{\m}$ and write
\beq
\label{equ:en1}
 P_{0}^{\m} \Bigl( \left\| p- p_{n}\right\| _{\infty } \ge {\delta}/{2}
\Bigr)     & \le &   n BM_{00} \Bigl( \sup_{0\leq  x \leq 1/n} \left|
p(x)\right|  \ge {\sqrt{\m} \delta }/{4} \Bigr)    \\ & \le &   2 n
BM_{0} \Bigl( \sup_{0\leq  x \leq 1/n} \left| p(x)\right| \geq  {\sqrt{\m}
\delta }/{4} \Bigr) \le    32 (\sqrt{n}/\d \sqrt{\m}) e^{ -  n \mu \d^2/
32}. \nonumber  
\eeq 
The first and third inequalities require no
explanation. For the middle inequality, note that if ${\cal A}$ is an event measurable
over  $\{ p(x), 0 \le x \le 3/4 \}$, then $BM_{00}( {\cal A} ) = 2 BM_0 [
e^{- p^2(3/4)/2}, {\cal A} ] \le 2 P_0 ( {\cal A} )$.

For the second term in (\ref{equ:enlarge}), bring in   $ {\cal I} (w)
=\frac{1}{2}\int_{0}^{1}$ $|w^{\prime}|^{2}$ the usual Brownian rate
function (${\cal I}(w) \equiv \infty $ when the integrand is not defined)
and  ${\cal I}(D) = {\inf }_{w \in D}  {\cal I}(w)$. Again, we first move
to the Bridge measure:   
\beqn 
P_{0}^{\m} \Bigl( \widehat{p}_{n} \in D
\Bigr)  & \leq  &   P_{0}^{\m} \Bigl( {\cal I}(\widehat{p}_{n}) \ge {\cal
I}(D) \Bigr) \\  & = & BM_{00} \Bigl( \frac{n}{2} \sum_{k=0}^{n-1}
|p((k+1)/n)-p({k/n})|^2 \ge  \m {\cal I}(D) \Bigr).   
\eeqn 
The latter
probability may be written out explicitly, and we continue to overestimate:
with $S_n(x) = (n/2)  ( x_1^2 + \sum_{k=1}^{n-2} (x_{k+1} - x_k)^2  +
x_{n-1}^2)$,  
\beq
\label{equ:en2}
\lefteqn{ BM_{00} \Bigl( \frac{n}{2} \sum_{k = 0}^{n-1}  |p({k+1}/n) -
  p({k/n})|^2 \ge \m {\cal I}(D) \Bigr) } \\ & = &  \frac{1}{(2 \pi/n
  )^{(n-1)/2}}   \int_{S_n(x) \ge \m {\cal I}(D)} \exp{\Bigl[ -S_n(x)
    \Bigr]} dx_1 \cdots dx_{n-1}   \nonumber \\    & \le & e^{-(1 - \eta)
  \m {\cal I}(D)}   \frac{1}{(2 \pi/n)^{(n-1)/2}} \int_{{\R}^{n-1}
} \exp{ \Bigl[- \eta S_n(x) \Bigr]} dx_1 \cdots dx_{n-1}  \nonumber \\  & =
&  \Bigl( \frac{1}{\sqrt{\eta}} \Bigr)^{n-1} e^{-(1 - \eta) \m {\cal
I}(D)}.  \nonumber  
\eeq 
Returning now to the event of interest
(\ref{equ:defD}), we combine  the two bounds  (\ref{equ:en1}) and
(\ref{equ:en2}) to find that 
\be
\label{equ:almost}
P_{0}^{\m}\Bigl(  D^{\mu}_{\d}(\gamma, \eta) \Bigr)  \le  32
 \frac{\sqrt{n}}{ \delta \sqrt{\m} } \exp \Bigl[ -n \mu \delta ^{2}/ 32
 \Bigr] +\left( \sqrt{\eta }\right)^{1-n} \exp \Bigl[ -(1-\eta ) \m   {\cal
 I}\Bigl( D_{\delta}^{\mu}(\gamma, \eta) \Bigr) \Bigr].  
\ee 
The first term
on the right can be made small by choosing $n$ appropriately.  It remains
to estimate $ {\cal I} ( D_{\delta}^{\mu})$ from below.   The results of
the Section $2.1$ imply that, in the present scaling,  
\be
\label{equ:varin}
   \frac{1}{2 \m} \int_0^1 {(g^{\prime})}^2  + \frac{1}{2}  \int_0^1 (1 -
       g^2)^2 \ge  \Bigl( \frac{8}{3} - \ep \Bigr) \frac{1}{\sqrt{\mu}} 
\ee
for any $g$ in $D_{\delta}^{\m}$ and $\ep = \epsilon(\mu) > 0$ going
to zero as $\m \uparrow \infty$.  This converts the problem into
one of bounding  $\int_0^1 (1-g^2)^2$ above.

Let $f \in D^{\m}$.  Directly from the definition of $D^{\mu}$ we see that
\[
\int f^4 - 2 \left(\int f^4\right)^{1/2}\! + 1 - 2 \mu^{-1/2} (\gamma+\eta)
\le 0,
\]
and the inequalities  
\be
\label{equ:upperlowerin}
1-\sqrt{2\left( \gamma +\eta \right)/\sqrt{\m }}\leq \left\| f\right\|
_{L^{4}}^{2}\leq 1+\sqrt{2\left( \gamma +\eta \right) / \sqrt{ \m }} 
\ee
follow immediately.  Next, if $g$ is to satisfy $\left\| f-g\right\|
_{\infty }<\delta $, then $\left| f+g\right| ^{2}<4\left|
f\right|^{2}+4\delta \left| f\right| +\delta ^{2}$ and so
\[
\int_{0}^{1}\left| f+g\right| ^{2} 
\le \left( \delta +2\left\| f\right\|_{L^{4}}\right) ^{2}
\]
by H\"older's inequality.  This last bound, together with
(\ref{equ:upperlowerin}), implies that
\[
\int_{0}^{1}\left| f+g\right|^{2} \le \left( \delta +2\sqrt{1+\sqrt{2\left(
\gamma +\eta \right)\ \sqrt{\m}}}\right) ^{2}<25
\]
for $\mu > 1$ and $\eta, \d, \gamma <1$, and so also 
\beqn
\int_{0}^{1}\left( 1-g^{2}\right) ^{2}-\int_{0}^{1}\left( 1-f^{2}\right)
^{2} & \leq & \delta \int_{0}^{1}\left| f+g\right| \left( \left|
1-f^{2}\right| +\left| 1-g^{2}\right| \right) \\  & \leq & 5\delta \left[
\sqrt{\int_{0}^{1}\left( 1-f^{2}\right) ^{2}}+\sqrt{
 \int_{0}^{1}\left(
1-g^{2}\right) ^{2}}\right].  
\eeqn after using Schwarz's
inequality. Therefore, $ || (1- g^2) ||_2 \le || (1- f^2) ||_2 + 5 \delta$
which, after squaring both sides and applying Cauchy's inequality, further
implies that \beqn \int_{0}^{1}\left(1 - g^{2} \right) ^{2}  & \le & \left(
1+\eta \right) \int_{0}^{1}\left( 1-f^{2}\right) ^{2} + 25 \left(
1+\eta^{-1} \right) \delta^{2} \\  & \leq & 2 \left( 1+\eta \right) \left(
\gamma +\eta \right) \frac{1}{\sqrt{\m}}  + 25 \left( 1+\eta ^{-1} \right)
\delta^{2} \eeqn for all $g \in D_{\d}^{\m}$.  Put together with
(\ref{equ:varin}), we have produced \be
\label{equ:Ilow}
{\cal I} \Bigl( D_{\delta }^{\mu}(\gamma, \eta ) \Bigr) \geq \left(
\frac{8}{3} - \ep -\left( \gamma +\eta \right) \left( 1+\eta \right)
\right) \sqrt{\mu} - \frac{25}{2} \left( 1+\eta^{-1}\right) \delta ^{2} \mu
\ee as the desired lower bound.

The final step revisits (\ref{equ:almost}) which, with the help of
(\ref{equ:Ilow}), says  \beq P_{0}^{\m} \Bigl( D^{\mu}(\gamma, \eta)
\Bigr) & \le & 32 \frac{\sqrt{n}}{\delta \sqrt{\m} } \exp{ \left[ \frac{-
\m n  \delta^{2}}{32} \right]} \\  &  & +\left( \frac{1}{\sqrt{\eta}}
\right)^{n}  \exp{ \Bigl[ \mu^2 \, \frac{25}{2} \left( 1+\eta^{-1} \right)
\delta^{2} \Bigr]} \exp{\Bigl[ -  \mu^{3/2} (1-\eta ) \left( \frac{8}{3} -
\ep  - ( \gamma +\eta)( 1+\eta) \right) \Bigr]}. \nonumber  \eeq A careful
choice of parameters  will now make all terms on the right negligible
compared to the one  of the form $\exp[{-\m^{3/2} etc.}]$. For example,
$\delta = {\m}^{-\alpha }$ and  $n \sim \m^{\beta }$ with $\frac{1}{4}$
$<\alpha <\frac{1}{2}$ and $\frac{3}{2}<\beta < \frac{1}{2 }-2\alpha $ (let
say $\alpha =5/16$ and $\beta = 5/4$) will do the job.  It follows that
\[
{\limsup_{\mu \ra \infty }} \, \mu ^{-3/2}\log P_{0}^{\m} \Bigl(
D^{\mu}\left( \gamma, \eta \right) \Bigr) \leq -(1-\eta )\left( \frac{8}{3}
- \left( \gamma +\eta \right) \left( 1+\eta \right) \right),
\]
and so, going back to the original quantity (see (\ref{equ:hd1}) through
(\ref{equ:defD})),  \beqn \limsup_{\mu \ra \infty}  \, \mu^{-3/2} \log
\int_{H^{\mu}\left( \gamma, \eta \right) }
e^{-\frac{1}{2}\int_{0}^{1}\left( \mu -p^{2}\right)^{2}} dP_{0}  & \leq &
-(1-\eta )\left( \frac{8}{3}-\left( \gamma +\eta \right) \left( 1+\eta
\right) \right) +\left( -\gamma +\eta \right) \\  & \leq & -(1-\eta )
\frac{8}{3} + \gamma \eta ^{2}+\eta \left( 2+\eta ^{2}\right) \eeqn for any
positive $\gamma < 1$.  At last, from the decomposition (\ref{equ:decomp}),
we deduce that
\[
\limsup_{\mu \ra \infty } \, \mu^{-3/2}\log \int_H e^{
\frac{1}{2}\int_{0}^{1}\left( \mu -p^{2}\right) ^{2}}dP_{0}\leq  -(1-\eta )
\frac{8}{3} + \frac{8}{3} \eta ^{2}+\eta \left( 2+\eta ^{2}\right)
\]
and letting $\eta \downarrow 0$ completes the proof.

\bigskip

\noindent{\bf Proof of Theorem {\ref{thm:rare}} }  We begin by following the
blueprint of the proof just completed. First, as in (\ref{equ:decomp}), the
integral over  $ H \backslash C_{\ep}^{\m}$ is first overestimated
by a sum of integrals according to the inclusion 
\[
   \Bigl( C_{\ep}^{\m} \Bigr)^c  \subseteq  \left\{ {\bigcup}_{0 \le k \le
(8/3 + \eta)/\zeta} \Bigl( H^{\m}( \zeta k, \zeta)  \cap (C_{\ep}^{\m} )^c
\Bigr)  \right\} \  \bigcup \ \left\{ p : \frac{1}{2} \int_0^1 (\m - p^2)^2
\ge (8/3 + \eta) \m^{3/2}  \right\}
\]
with any $\zeta > 0$.  The integral of $\exp[-(1/2) \int_0^1 (\m - p^2)^2 ]$ over
the last set on the right hand side trivially satisfies the desired bound.
Next, as in (\ref{equ:hd}) and (\ref{equ:hd1}), 
bounding the other integrals in this decomposition
comes down to bounding the $P_0$ probabilities of the
events $\{ H^{\m}( \zeta k, \zeta)  \cap (C_{\ep}^{\m} )^c \}$.  Following the
proof of Theorem \ref{thm:LD} further we come to the critical point. By
comparison with (\ref{equ:varin}) and the surrounding discussion we see that we 
now need an improved version of that variational inequality. In particular,
we require a lower bound of the form 
\[
      \frac{1}{2 \m} \int_0^1 (g^{\pr})^2 + \frac{1}{2} \int_0^1 (1 -
                 g^2)^2 \ge (8/3 + {\tilde{\eta}})  \frac{1}{\sqrt{\m}}
\]  
where
${\tilde{\eta}}= {\tilde{\eta}}(\ep) > 0$ depends on $\ep$ but is
fixed for all $g$ restricted to lie in $(C_{\ep}^{\mu})^c$,
appropriately scaled and $\delta$-enlarged.  
With the appropriate scaling, it is equivalent to show that
\be
\label{equ:scalelower} 
  \liminf_{a \ra \infty}  \left( {\inf} \Bigl\{  I(f; a),  \ f \in H_a \cap
    d(f, \{ f_a^* \}) > \ep \Bigr\} \right) > 8/3 
\ee 
where $I(f;a)$ and $H_a$ are as defined in (\ref{equ:scaledvar})  and 
$\{ f_a^* \}$ represents the set of minimizers of $I(f;a)$. 
    
We argue by contradiction. If (\ref{equ:scalelower})  failed to hold,
we could find a sequence $\{ {\tilde{f}}_a \}$  satisfying the
constraints, but so that $ I( {\tilde{f}}_a; a) \ra 8/3$.  By virtue of
the fact that $\int_{-a}^{a} {\tilde f}_a =  0$,  each $\tilde{f}_a$ 
has at least two zeros, $z_a^1 < 
z_a^2$. Further,  $|z_a^1 - z_a^2| \ra \infty$ and $(2a - |z_a^1 -
z_a^2| ) \ra \infty$.  If instead
 all the zeros were contained in a fixed interval $I$ 
(which may be assumed to be centered about the origin), it would follow
that  $|\int_{I^c}
{\tilde{f}}_a |$ must exceed a positive  multiple of $a$.  Indeed,
$\int_{-a}^{a} [1- ({\tilde{f}}_a)^2]^2$ remains bounded and the length
of $I^c$ is itself $O(a)$. But then  
$|\int_I {\tilde{f}}_a| \ge const. \times  a$ to maintain the mean zero 
condition causing $|{\tilde{f}}_a| \ge const. \times a$ on some subset 
of $I$ of positive measure.  This in turn would imply that $\int_I [1
- ({\tilde{f}}_a)^2]^2$ 
grows without bound as $a
\ra \infty$, and that is impossible.

At this point we return to the strategy behind the proof of Theorem
{\ref{thm:rate}}.  First we fix large symmetric intervals of length $2M$ 
around each of the two zeros $z_a^1$ and $z_a^2$ specified thus far.
When considering these intervals separately we will identify 
$z_a^1$ or $z_a^2$ with the origin as may be done by translation. 
By the core argument behind the proof of Theorem {\ref{thm:rate}}, on each of 
these intervals we can find subsequences $\{ \tilde{f}_a^1 \} $ and 
$\{ \tilde{f}_a^1 \} $ converging to $ f_{\infty}^1$ and $f_{\infty}^2$
respectively.  Again, both   $f_{\infty}^1$ and
$f_{\infty}^2$ lie in $ H^1 \cap L^{\infty}$. Also, since both converge 
weakly in $H^1$, by lower semi-continuity of the functional $I$, it follows
that $I(f_{\infty}^{k}, \infty) \ge 4/3$, and both also satisfy the equation
$f_{\infty}^{\prime} = \pm (1 - f_{\infty}^2)$  on their respective
domains.  That is,  $f_{\infty}^k(x) = \pm \tanh(x)$ for $k =1,2$.

Next we show that $z_a^1$ and $z_a^2$ are in fact the only zeros of
${\tilde f}_a$ for $a \ra \infty$, and that $z_a^2 - z_a^1$ is
roughly $a$:
\be
\label{equ:zeros}
\lim_{a \uparrow \infty} \Bigl( |z_a^1 - z_a^2| - {a} \Bigr) = 0.
\ee
For the first part, simply note that if there were a third zero we
may repeat the argument of the preceding paragraph to conclude that the rate function $I$
would exceed $4/3 + 4/3 +4/3$ as $a \ra \infty$, but that contradictions
the assumption $I({\tilde f}_a, a) \ra 8/3$.  
For the statement regarding the distance between $z_a^1$ and $z_a^2$, 
note first that the integrals of  $\{\tilde{f}_a\}
$ over $ [ ( z_a^1 + z_a^2)/2 ,  ( z_a^1 + z_a^2)/2 - a ]$ or over its
complement in $[-a, a)$ must have the same absolute value but opposite
signs. Translating to place $z_1^a$ at the origin  we find that
\[
\int^{u}_{\ell} {\tilde f}_a(x) dx  + \int^{-\ell}_{-u} {\tilde f}_a(x) dx
= 0
\]
where  $u = ( z_a^2 -
z_a^1)/2$ and $\ell = ( z_a^2 - z_a^1)/2 - a$. We can now assume that $f_{\infty}^1
= \tanh$ so that $f_{\infty}^2 = -\tanh$. Then, since there is uniform
convergence to these limiting functions, both integrals in the last 
display are approximately $u + \ell$. It follows that
 $ u + \ell \ra 0$ as $a \rightarrow \infty$ which is the same as (\ref{equ:zeros}).

The conclusion is that $||{\tilde f}_a - f_a||_{\infty} \ra 0$
where $f_a$ is the test function constructed in (\ref{equ:test}).
Since $f_a$ minimizes $I(f;a)$ as $a \ra \infty$ we have shown that 
$d({\tilde f}_a, \{ f_a^* \} )
\ra 0$,  contradicting the original hypothesis. The proof is
complete.

\section{Expanding about the extrema}
\setcounter{num}{3} \setcounter{equation}{0}
\label{sec:expand}

Following the classical Laplace method we wish to expand $f(\m)$ in the
vicinity of each $I_{\m}$-minimizer.  While we have not actually computed
any such minimizer at finite $\m$ (nor have we proved the anticipated
uniqueness up to translation),  it suffices to introduce the following
proxy.  The Euler-Lagrange equation (\ref{equ:EL1}), describing the scaled
minimizer(s), may be solved in terms of Jacobi elliptic functions.  Thus
motivated, we bring in \be
\label{equ:pmu}
p_{\m}(x) = k \sqrt{\m}  \times \sn( \sqrt{\m} x, k)  \ee along with its
translates $p_{\m}^a(\cdot) = p_{\m}(\cdot + a)$.  As a function of the
real variable $x$, $\sn(x,k)$ is periodic with period determined by its {\em
modulus} $k \in [0,1]$.  In particular, $\sn(\cdot, k) = \sn(\cdot + 4K, k)$ in
which $K$ is the complete elliptic integral of the first kind: $ K(k) =
\int_0^1 [(1-x^2)(1-k^2x^2)]^{-1/2} dx$.  (For background on {\em sin-amp}
and other elliptic functions used throughout,  \cite{BF} is recommended.)

In $p_{\m}$, we must choose $k$ so that $4 K = \sqrt{\mu}$, and it may be
deduced that  $ k^2 \simeq 1 - 16 e^{-\sqrt{\m}/2}$.  With this parameter
set, an exact computation will yield \be
\label{equ:rate__p}
   I_{\m} ( p_{\mu}) = \frac{8}{3} \m^{3/2} + O \Bigl( e^{-\sqrt{\m}/4}
\Bigr).  \ee While this is certainly heart-warming, a more important
connection with the discussion in Section 2.1 is seen in the fact that
$\sn(x,k) \simeq \tanh(x)$ for $ k \up 1$. So, with $p_{\m}^*$ any
$I_{\m}$-minimizer with $p_{\m}^*(0) = 0$, the proofs of Theorem 2.1
and 2.3 will explain why
\[
\label{pmp}
   \lim_{\m \ra \infty} \Bigl| \Bigl| \frac{1}{\sqrt{\mu}} p_{\m}  -
                                   \frac{1}{\sqrt{\m}} p_{{\m}}^*  \Bigr|
                                   \Bigr|_{\infty}  = 0.
\]
In other
words, an appropriate $L^{\infty}$  tube about the set of translates $\{
p_{\m}^a \}$ contains a like tube about the set of $I_{\m}$-minimizers.
Theorem \ref{thm:rare} then implies the following.

\begin{cor}
\label{cor:newcenter}
Let $\ep > 0$. Then \be
\label{equ:errter}
  f(\m) = \frac{1}{\sqrt{2 \pi}} E_0 \Bigl[ e^{-\frac{1}{2} \int_0^1 (\m -
p^2)^2} A(p) , \ d( p, \{p_{\m}^a \}) \le \ep \sqrt{\mu} \Bigr] + O \Bigl(
e^{-(8/3+\eta) \m^{3/2}} \Bigr), \ee with some  $\eta = \eta(\ep) > 0$.
\end{cor}

This last observation turns the problem of expanding $f(\mu)$ 
about each $I_{\m}$-minimizer into that of expanding the
expectation in (\ref{equ:errter}) about each $p_{\m}^a$, $0 \le a < 1$.
While a definite advancement, we must still confront
the  degeneracy inherent in the translation
invariance of $I_{\m}$.  We handle this issue a conditioning
procedure in order to pin the $E_0$ expectation about a single $p_{\m}^a$,
which for convenience  is taken to be $p_{\m} = p_{\m}^0$. The idea is
that since the translations  $p_{\m}^{0} \ra p_{\m}^{0+ \ep}$ are generated
by $p_{\m}^{\pr} =  const. \times \cn(\sqrt{\m} x , k) \dn(\sqrt{\m} x,
k)$,  conditioning the path to be orthogonal to $p_{\m}^{\pr}$ will keep
the path in a small neighborhood of $\pm p_{\m}$ as opposed to 
translates further afield.  Relating the $E_0$ expectation to this
conditioned version of itself requires a change of measure 
formula provided in the following lemma; the proof is deferred to the appendix.

\begin{lemma}
\label{l:StatZeros}
Suppose $X(\cdot)$ is a smooth stationary process, periodic, of period one.
Also assume that, with probability one, $X$ has at least one zero. If ${F}$
is a functional  that is invariant under translations of the path, then
\[
  E \Bigl[ { F}(X) \Bigr] =  E \Bigl[ { F}(X) \frac{1}{{\cal N}^{-1}}
  \Bigl|  X(0)=0 \Bigr]  P \Bigl( X(0) = 0 \Bigr)
\]
where
\[
{\cal N} = {\cal N}(X) =  \sum_{z \in {\cal Z}} |X_z^{\pr}|^{-1}
\]
and ${\cal Z}$ is the set of zeros of $X$.  Here and throughout,
the notation $P(X = a)$ indicates the  density of $X$ at $a$.
\end{lemma}

Applying Lemma {\ref{l:StatZeros}} to the matter at hand, we 
set\footnote{The choice of notation will become clear in the next section.}
\[
   \phi_1^{\m}(x) \equiv \frac{ \cn( \sqrt{\m} x, k) \dn( \sqrt{\m} x, k) }
                              { \sqrt{ \int_0^1   \cn( \sqrt{\m} x^{\pr},
                              k)  \dn( \sqrt{\m} x^{\pr}, k) dx^{\pr}}},
\] 
and note the following.

\begin{cor} 
\label{c:shift1}
Let 
\be
\label{equ:defR}
R(p) \equiv {\cal N}^{-1} \left( x \ra \int_0^1 \phi_1^{\m}(x +
                                x^{\pr}) p(x^{\pr})  d x^{\pr} \right),
                                \ee then \beqn \lefteqn{
\hspace{-1cm} E_0 \Bigl[   e^{-\frac{1}{2} \int_0^1 (\m - p^2)^2} A(p),
                  d(p, \{ p_{\m}^a \}) \le \ep \sqrt{\m} \Bigr]}  \\ & = &
                  E_0^0 \Bigl[   e^{-\frac{1}{2} \int_0^1 (\m - p^2)^2}
                  A(p) R(p) ,   d(p, \{ p_{\m}^a \}) \le \ep \sqrt{\m}
                  \Bigr] \ P_0 \Bigl( \int_0^1  \phi_1^{\m} p = 0 \Bigr).
\eeqn 
Here $E_0^0$ is now the $CBM$ conditioned so that
both $\int_0^1 p = 0 $ and  $\int_0^1   \phi_1^{\m} p = 0$.
\end{cor}

Next, consider the intersection of $\{ p : d(p, \{p_{\m}^a\}) \le \ep
\sqrt{\m} \}$ and  $\{ p : \int_0^1 \phi_1^{\m} p = 0 \}$.  It is easy to
see that the resulting set will contain the union of $\{ ||p -
p_{\m}^0||_{\infty } \le \ep_1 \sqrt{\m} \}$  and the same object with
$p_{\m}^{0}$ replaced by  $p_{\m}^{1/2} = - p_{\m}^{0}$ for some $\ep_1$
small enough.  Likewise, it will be contained in a similar union with
$\ep_1$ replaced by some larger $\ep_2 >0$.  Of course, the integral in
question is invariant under the sign change $p \ra -p$.  These comments
along with  Corollaries \ref{c:shift1} and \ref{cor:newcenter} allow the
following statement: 
\be
\label{equ:apxf}
  f(\mu) \simeq  \sqrt{ \frac{2}{\pi}}  E_0^0 \Biggl[  e^{-\frac{1}{2}
  \int_0^1 (\m - p^2)^2} A(p) R(p), || p - p_{\m} ||_{\infty} \le \ep
  \sqrt{\m}  \Biggr]  P_0  \Bigl( \int_0^1  \phi_1^{\m} p = 0 \Bigr) 
\ee 
up to $O(e^{-(8/3+) \m^{3/2}})$ errors granted that we eventually  prove
that the desired level of asymptotics for the $E_0^0$ integral in
(\ref{equ:apxf})  are independent of $\ep > 0$.

Having centered the integral about the single path $p_{\m}$, we  complete
this section by performing the change of variables  $p \ra p + p_{\m}$ in
order to bring the contribution of  $p_{\m}$ up into the exponent.

\begin{prop}  
\label{p:shift}
With the $E_0^0$ integral on the right hand side of (\ref{equ:apxf})
denoted by $f(\m; \ep)$ we have \be
\label{equ:shift}
   f(\m ; \ep) = e^{- I_{\m}(p_{\m}) } E_0^0 \Biggl[  e^{-\frac{1}{2}
                      \int_0^1 ( q_{\m} - 2 \m) p^2} e^{-2 p_{\m} \int_0^1
                      p_{\m} p^3 - \frac{1}{2} \int_0^1 p^4} A(p + p_{\m})
                      R(p + p_{\m}), || p ||_{\infty} \le \ep \sqrt{\mu}
                      \Biggr] \ee in which $q_{\m}(x) \equiv 6 \m  k^2
                      sn^2(\sqrt{\m} x, k)$ $( = 6 p_{\m}^2(x))$.
\end{prop}

We have thus extracted the advertised leading term, $e^{-I_{\m}(p_{\m})}
= e^{-8/3 \m^{3/2}}(1+ o(1))$.
Further, the formula (\ref{equ:shift}) identifies  the Gaussian measure
$e^{-(1/2) \int_0^1 (q_{\m} - 2\m) p^2} d P_0^0$ which,  as is the case
in finite
dimensional Laplace asymptotics, will dictate the remainder of our
computation.   
The study of this measure is initiated in the next section.

\bigskip

\noindent{\bf Remark}  Given (\ref{equ:shift}), it is a simple
matter to obtain the lower bound complementing Theorem \ref{thm:LD}: $
\lim_{\m \ra \infty} {\m}^{-3/2} \log f(\m) = - 8/3$.

\bigskip

\noindent{\bf Remark} If we now understand that a vicinity of $p_{\m}$ (and
its translates) accounts  for the leading order behavior of $f(\m)$ for  $\m \ra
\infty$, it is interesting to consider what this implies for  the random
potential. Running the Ricatti correspondence ``backwards'' relates this  leading
path to the potential $q(x; \mu) \equiv - \m + p_{\m}^{\pr}(x) + p_{\m}^2(x) $
$\simeq - 2 \m \, \sech^2( \sqrt{\m}(x - 1/2) )$. While formal, this
indicates that large negative deviations of the ground state stem from
White-Noise potentials lying nearby a single well of depth $\m$ and width
$1/\sqrt{\m}$.

\bigskip

\noindent{\bf Proof of Proposition {\ref{p:shift}} }  This is a consequence of the 
Cameron-Martin formula for $P_0^0$
proved in the Appendix (Lemma \ref{l:CMart}).  It states that, for bounded functions
$F$ of the path,
\[
 E_0^0 \left[ F{\left( p\right)} \right] =
     E_0^0  \left[ F{\left( p+p_{\m} \right)}
\exp \left\{ \int_{0}^{1}{p_{\m}^{\pr \pr }} p-\frac{1}{2}
\int_{0}^{1}\left| {p_{\m}^{\prime}} \right|^{2} \right\} \right].
\]
In the present case 
\[
F(p) = e^{-\frac{1}{2} \int_0^1 (\m - p^2)^2} A(p) R(p) 1_{\{ || p -
  p_{\m}||_{\infty} \le \ep \sqrt{\m}\}},
\]
and a simple expansion yields
\beqn
\lefteqn{ \frac{1}{2} \int_0^1 \Bigl( \m - (p + p_{\m})^2 \Bigr)^2   +
  \frac{1}{2} \int_0^1 |p_{\m}^{\pr}|^2  
  - \int_0^1 p_{\m}^{\pr \pr} p 
             } \\
& = &     I_{\m}(p_{\m})  + 2 \int_0^1 p_{\m} (p_{\m}^2 - \m) p  
          + \frac{1}{2} \int_0^1 \Bigl(6 p_{\m}^2  - 2 \m +  4 p_{\m} p +
	  p^2 \Bigr) p^2  - \int_0^1 p_{\m}^{\pr \pr} p.
\eeqn
Next, we notice that $ p_{\m}^{\pr \pr} = 2 p_{\m}^3 - 2 \m p_{\m} + constant$, and therefore
\[
  - \int_0^1 p_{\m}^{\pr \pr} p = - 2 \int_0^1  p_{\m} ( p_{\m}^2 - \m) p   
\]
when $\int_0^1 p = 0$.
The proof is finished by combining the last two formulas.

\section{Hill's Spectrum and the associated Gaussian process}
\setcounter{num}{4}
\setcounter{equation}{0}
\label{sec:hillspec}

Our analysis of the culminating form of the expectation (\ref{equ:shift}) takes the anticipated route.
Restricted to a set of relatively small $L^{\infty}$ norm, it is expected that  $A(\cdot)$ and $R(\cdot)$
settle down to $A(p_{\m})$ and $R(p_{\m})$ as $\mu \ra \infty$.  More delicate, the cubic $(2 \int_0^1 p_{\m} p^3)$ 
and quartic $( -(1/2) \int_0^1 p^4)$ terms in the exponential should be 
lower order when compared with the quadratic factor $ (1/2) \int_0^1 ( q_{\m} - 2 \m) p^2$ 
which is of order $\m$.  That is,  
the remainder of the computation should be viewed with 
respect to the Gaussian
measure $ \exp[ -(1/2)  \int_0^1 ( q_{\m} - 2 \m) p^2 ] d P_0^0$
arising from the Hessian of the rate function $I_{\m}(p)$ at $p = p_{\m}$.

Of course, exercising this program first requires being able to compute in the latter measure.   This prompts
the study of the periodic spectrum of the (deterministic) operator
\[
   Q_{\m} =  - \frac{d^2}{dx^2} + q_{\m}(x), 
\] 
and in this we are met with no small piece of good fortune.  The reason we can compute the details
of $f(\m)$ rests on the rather special properties of $Q_{\m}$.

\subsection{Hill's equations}

Consider the general Hill's operator
$Q = - d^2/d x^2 + q(x)$
in which $q(x)$ is smooth and with period 
now taken to be $1/2$ (note our motivating example $Q_{\m}$).
The periodic spectral points of $Q$ comprise a list:
\[
   - \infty < \ld_0 < \ld_1 \le \ld_2 < \ld_3 \le \ld_4 < \ld_5 \le \ld_6 < \ld_7 \le \cdots \uparrow + \infty.
\] 
In particular, the so-called principal series $\ld_0 < \ld_3 \le \ld_4 < \ld_7 \le \ld_8 < \cdots$ 
makes up the periodic 
spectrum of $Q$ acting on $L^2$ functions of period $1/2$, and the complementary series 
$\ld_1 \le \ld_2 < \ld_5 \le \ld_6 < \cdots $ fills out the periodic spectrum of $Q$ on $L^2$ functions
of period $1$. 
An equivalent
characterization of spectrum may be described with the help of Hill's discriminant
\[
   \Delta(\ld) = y_1(1/2,\ld) + y_{2}^{\pr}(1/2,\ld)
\] 
in which $y_1(1/2,\ld)$ ($y_2(1/2,\ld)$) is the normalized sine-like (cosine-like) solution 
of $Qy = \ld y$ with $y(0) = 0, y^{\pr}(0) = 1$ ($y(0) = 1$, $y^{\pr}(0) = 0$).  The classical result is
that $\Delta(\ld)$
is an entire function of order $1/2$ and that it encodes the spectrum: 
$\Delta(\ld) = + 2$ on the principal series and 
$\Delta(\ld) = -2$ on the complementary series.

A special situation occurs when the shape of the potential $q$ is such that
the simple eigenvalues $\ld_0 < \ld_1 < \cdots < \ld_{2g}$ are finite in number 
with the rest of the list double: $\ld_{2\l-1} = \ld_{ 2 \l}$ for $ \l > g$.   
Then $Q$ is said to be {\em finite gap}, and it is the remarkable
discovery of Hochstadt \cite{Hoch} that in this case  
the simple spectrum determines the full spectrum and so also $\Delta(\ld)$.
More precisely:

\vspace{.2cm}

\noindent{\bf Hochstadt's Formula} {\em Let} $Q$ {\em be finite-gap 
with} $2 g + 1$ {\em simple eigenvalues}.  {\em Then} $\Delta(\ld) = 2 \cos \psi(\ld)$  {\em with}
\be
\label{equ:hoch}
   \psi(\ld) = \frac{\sqrt{-1}}{2} \int_{\ld_0}^{\ld}  
\frac{(s - \ld_1^{\pr})\cdots(s - \ld_g^{\pr})}{  \sqrt{ - (s -\ld_0) \cdots (s -\ld_{2g})}} d s
\ee
{\em in which} $\ld_1^{\pr} < \cdots < \ld_g^{\pr}$  {\em are the points}
$\ld_{2\l-1} < \ld_{\l}^{\pr} < \ld_{2\l}$ {\em where} $\Delta^{\pr}(\ld) = 0$.  {\em They are determined 
  from the simple spectrum through the requirement:} $\psi(\ld_{2\l}) - \psi(\ld_{2\l-1}) = 0$ {\em for} 
  $\l = 1,2, \dots g$. 

\vspace{.2cm}

This formula will play an essential role in Sections 6 and 7.  More background information
on Hill's spectrum can be found in \cite{MW} or \cite{MvM}.

\subsection{When $Q = Q_{\m}$}

The key fact is that the family of Lam\'e operators 
$Q = - d^2/dx^2 +  m(m+1) k^2 \sn^2( x, k)$ over the period $0 \le x \le 2 K$ 
are finite gap with $g = m$ (see again \cite{MvM}).  In $Q_{\m}$ we have $m=2$, and, what is
more,  the simple eigenvalues and
corresponding eigenfunctions are known, having first 
been computed in \cite{I}.  With 
\[ 
a_{\pm}(k) = 1 + k^2 \pm \sqrt{1-k^2 +k^4},
\]
we have:  
\beq
\label{equ:eigs}
    \ld_0 = 2 a_{-}(k) & & \tilde{\phi_0}(x) = 1 -  a_{-}(k) \sn^2( x , k)  \\
   \ld_1 =   1 + k^2  & & \tilde{\phi}_1(x) = \cn(x, k) \dn(x, k) \nonumber \\
   \ld_2 =  1 + 4 k^2  & &   \tilde{\phi}_2(x) =  \sn(x, k) \dn(x, k)  \nonumber \\
    \ld_3 =   4 + k^2  & &  \tilde{\phi}_3(x) = \sn(x, k) \cn(x, k) \nonumber \\
   \ld_4 =  2 a_{+}(k)  & &
   \tilde{\phi}_4(x) = 1 - a_{+}(k) \sn^2(x, k).  \nonumber
\eeq
Of course, things are to be scaled as in $ x \ra 4K x = \sqrt{\mu} x$ 
to keep the period at $1/2$.  Accordingly,
$\ld_{\l} \ra \ld_{\l}^{\m} = \m \times \ld_{\l}$  and  
$\tilde{\phi}_{\l}(x) \ra \tilde{\phi}_{\l}^{\m}(x) = \tilde{\phi}_{\l}(\sqrt{\m} x)$,
after which we may introduce the 
$L^2[0,1]$-orthonormal sequence 
denoted by $ \{ \phi_{\l}^{\m} \} =$ 
$ \{ \tilde{\phi}_{\l}^{\m} /  || \tilde{\phi}_{\l}^{\m} ||_{2}  \}$.
Keep in mind though the unscaled $\ld$'s and $\tilde{\phi}$'s depend on $\m$
through the modulus $k$.   Further, it is worthwhile 
noting that unscaled operator $\m^{-1/2} Q_{\m}(\cdot/\sqrt{\m})$
tends to $-d^2/dx^2 + 6 \tanh^2(x)$ over the whole line as $\m \ra \infty$ or $k \ra 1$.  
In this picture,
$\tilde{\phi}_0, \dots, \tilde{\phi}_3$ correspond to bound states 
($\tilde{\phi}_0, \tilde{\phi}_1 \simeq \sech^2$
and $\tilde{\phi}_2, \tilde{\phi}_3 \simeq \sinh \sech^2$
up to $O(e^{-\sqrt{\m}/2})$ errors in $L^{\infty}$)  with continuous spectrum beginning at 
$\ld_4 \simeq 6$.

\subsection{The Gaussian Measure}

We now have a more complete picture 
connecting the basic degeneracy and the introduced conditioning.
If the $CBM$ had appeared unconditional in the definition of $f(\mu)$, we would
at this point be faced with the measure $\exp[ - 1/2 \int_0^1 (q_{\m} - 2 \m) p^2 ] \times d CBM$
to provide concentration of the path about the extrema set.
However, the first two members of the corresponding  
spectrum  $\ld_0^{\m} - 2 \m$ and $\ld_1^{\m} - 2 \m$ are $ \simeq - 16 e^{-\sqrt{\m}/2}$, 
preventing this measure from having a sense.\footnote{That these first two 
spectral points  are even a little negative
is due to the fact that $p_{\m}$ is not actually an extrema, 
though exponentially close.}  
Thus, the mean-zero conditioning built into the problem works
to counteract the first
degeneracy tied to the ground state of $Q_{\m}$.  Furthermore, the 
conditioning we introduced to account for the translational degeneracy exactly 
removes the second eigenstate of $Q_{\m}$. 

That the imposed  $\int_0^1 p = 0$  only ``works to counteract''
the degeneracy at $\phi_0^{\m}$ is because that mode is not fully 
removed by the conditioning.  The constant function can however
be written in the eigenbasis: from 
the first and last items of the list (\ref{equ:eigs}) we see that
\be
\label{equ:one}
    1 = c_0  \phi_0^{\mu} (x)  - c_4 \phi_4^{\mu} (x)
\ee
where 
\be
\label{equ:c's}
  c_0 =  c_0(\m) = \Bigl(\frac{a_{+}(k)}{a_{+}(k) - a_{-}(k)} \Bigr)  || \tilde{\phi}_0^{\m} ||_2 
   \ \ \mbox{ and } \ \  c_4 = c_4(\m) = \Bigl( \frac{ a_{-}(k)}{ a_{+}(k) - a_{-}(k)} \Bigr) || 
         \tilde{\phi}_4^{\m} ||_2.
\ee
This allows an explicit description of the bottom of the 
conditioned spectrum. 

\begin{lemma} 
\label{lem:starmeas}
The Gaussian measure $P_{\m}^*$ defined by 
\[
   E_{\m}^* \Bigl[ F(p) \Bigr] = \frac{ E_0^0 
\Bigl[ F(p) \exp{ \Bigl\{ - \frac{1}{2} \int_0^1 (q_{\m}(x) - 2 \m) p^2(x) dx \Bigr\} }  \Bigr ]}
 {E_0^0 \Bigl[  \exp{ \Bigl\{  - \frac{1}{2} \int_0^1 ( q_{\m}(x) - 2 \m ) p^2(x) dx \Bigr\}  }   \Bigr]}  
\]
has the expansion
\beq
\label{equ:pexp}
p(x)  & = &   \Bigl( 
\frac{ c_4 \phi_0^{\m}(x) + c_0 \phi_4^{\m}(x) }{ \sqrt{ c_4^2 ( \ld_0^{\m} - 2 \m ) 
                           + c_0^2 (\ld_4^{\m} - 2 \m) }}  \g_0  + 
                            \frac{\phi_2^{\m}(x)}{\sqrt{\ld_2^{\m} - 2 \m}} \g_2
                 + 
                            \frac{\phi_3^{\m}(x)}{\sqrt{\ld_3^{\m} - 2 \m}} \g_3  \Bigr) \\
&  & 
                           +  \Bigl( \sum_{\l \ge 5}  
                          \frac{1}{\sqrt{{\ld}_{\l}^{\m} - 2 \m   }}   
                            \phi_{\l}^{\m}(x)  \g_{\l} \Bigr) \nonumber \\
& \equiv &   p_{\l}(x) \, + \, p_h(x) \nonumber
\eeq
for $\{ \g_{\l} \}$  a sequence of independent standard Gaussians.  
\end{lemma}

The division into low and high modes ($p_{\l}$ and 
$p_{h}$) is prepared for later.  We pause here to note that 
while the typical normalizer of each high mode $\sqrt{\ld_{\l}^{\m} - 2 \m}$
is $ O(\sqrt{\m})$ (for any $\l \ge 2$), the corresponding object in
the lowest mode $\sqrt{   c_0^2 (\ld_4^{\m} - 2 \m) + etc.}$ is only
of order $\m^{1/4}$.  Thus, 
scaling out the $\sqrt{\m}$ from the exponential weight of $P_{\m}^*$
(as is customary in Laplace type computations), we see that the
resulting Gaussian measure has a spectral gap disappearing like $\m^{-1/4}$.   
This is degeneracy (or perhaps it is better to say ``near-degeneracy'') 
orthogonal to the set of extrema mentioned in the introduction.  
It would be interesting to understand its physical significance.

\bigskip

\noindent{\bf Proof }  
This is just the classical Karhunen-Lo$\grave{\mbox{e}}$ve expansion,
see \cite{A}; we outline the derivation to make clear what  occurs
at the bottom of the spectrum on account of the condition $\int_0^1 p = 0$.

As $CBM[ \exp{\{- 1/2 \int_0^1 (q_{\m} - 2\m) p^2} \}] = +\infty$,
it is convenient to express things through $ d P_{q_{\m}} = 
Z_{q_{\m}}^{-1} \exp{[- 1/2 \int_0^1 q_{\m} p^2]} d CBM$
which does have finite total mass. That is, we write
\be
\label{equ:renorm1}
 E_{\m}^* \Bigl[ F(p) \Bigr] = \frac{ E_{q_{\m}}  
\Bigl[ F(p) \, e^{ \m  \int_0^1  p^2(x) dx  }  \, \Bigl| \, \int_0^1 p(x) = 0, \int_0^1 p(x) \phi_1^{\m}(x) = 0 \Bigr]}
 {E_{q_{\m}} \Bigl[  e^{ \m  \int_0^1  p^2(x) dx } \, \Bigl| \, \int_0^1 p(x) = 0, \int_0^1 p(x) \phi_1^{\m}(x) = 0   \Bigr]},  
\ee
and introduce the coordinates $ p(x) = \sum_{\l = 0}^{\infty} \frac{1}{\sqrt{\ld_{\l}^{\m}}} \phi_{\l}^{\m}(x) \eta_{\ell}$
under $P_{q_{\m}}$ (the $\eta_{\l}$'s  are independent standard Gaussians).  
Written in this way,
the numerator of (\ref{equ:renorm1}) takes the form
\beqn
\lefteqn{
  E \Bigl[ F \Bigl( \sum_{\l = 0}^{\infty} \frac{1}{\sqrt{\ld_{\l}^{\m}}} \phi_{\l}^{\m}(x) \eta_{\ell}
   \Bigr) \, \exp{ \Bigl\{ \sum_{\l = 0}^{\infty} \frac{\m}{\ld_{\l}^{\m}}  \eta_{\l}^2  \Bigr\} }  
   \, \Bigl| \,  \frac{c_0}{\sqrt{\ld_0^{\m}}}  {\eta}_0 
           -  \frac{c_4}{\sqrt{\ld_4^{\m}}} {\eta}_4 = 0, \, \eta_1  = 0 \Bigr]} \\
& =   & C \, E \Bigl[  F \Bigl( \frac{1}{\sqrt{\ld_0^{\m}}} 
                             \phi_0^{\m}(x) \eta_0 + 
\frac{1}{\sqrt{\ld_4^{\m}}} \phi_4^{\m}(x) \eta_4 +  {\tilde p}_h(x)
\Bigr) 
         \exp{ \Bigl\{\frac{\m}{\ld_{4}^{\m}}  \eta_{1}^2  + \frac{\m}{\ld_{4}^{\m}}  
                             \eta_{4}^2 \Bigr\}} 
                               \Bigl| \,  \frac{c_0}{\sqrt{\ld_0^{\m}}}  {\eta}_0 
           -  \frac{c_4}{\sqrt{\ld_4^{\m}}} {\eta}_4 = 0    \Bigr].
\eeqn
Here $C$ is a constant factor and
\[
     {\tilde p}_h(x) =   \frac{\phi_2^{\m}(x)}{\sqrt{\ld_2^{\m} - 2 \m}} \g_2 + 
                 \frac{\phi_3^{\m}(x)}{\sqrt{ \ld_3^{\m} - 2 \m}} \g_3 + p_h(x),
\]
its distribution now identified by independence.

There remains the distribution of 
$p_{0,4}(x) =  \frac{1}{\sqrt{\ld_0^{\m}}} \phi_0^{\m}(x) \eta_0 + 
     \frac{1}{\sqrt{\ld_4^{\m}}} \phi_4^{\m}(x) \eta_4 $
subject to the quadratic weight $\exp[ \frac{\m}{\ld_{4}^{\m}}  \eta_{1}^2  + \frac{\m}{\ld_{4}^{\m}}  
                             \eta_{4}^2]$ and 
linear conditioning $ \{  \frac{c_0}{\sqrt{\ld_0^{\m}}}  {\eta}_0 
           -  \frac{c_4}{\sqrt{\ld_4^{\m}}} {\eta}_4 = 0 \} $. 
The question is only whether
the latter conditioning can hold down the focusing weight.  A straightforward computation
shows that, for $X, Y \sim N(0,1)$ and independent,  
\[
 (1 - 2 \a ) + \Bigl(\frac{a}{b}\Bigr)^2 (1 - 2 \b) > 0  \; \; \mbox{ implies } \; \;   
E[ e^{ \a X^2 + \b Y^2} \, | \, a X + b Y = 0] < \infty.
\]
Next note that, by inspection, 
$a_{+}(k) \simeq 3$ and $ a_{-}(k) \simeq 1$ 
up to errors of order $1-k^2 = O(e^{-\sqrt{\m}/2})$. 
At the same level approximation we also have
\[
  || \tilde{\phi}_0 ||_2^2  \simeq \frac{1}{\sqrt{\m}}
  \int_0^{\sqrt{\m}} \cn^4(x, k) dx  \simeq \frac{2}{\sqrt{\mu}} \int_{-\infty}^{\infty} \sech^4(x) dx 
                                   =  \frac{8}{3 \sqrt{\m}}, 
\]
and
$ || \tilde{\phi}_4 ||_2^2 \simeq  4;$
the upshot being that $c_0(\m) \simeq \sqrt{6} \m^{-1/4}$ and $c_4(\mu) \simeq 1$. Thus,
\[
  \Bigl(1 - \frac{2 \m}{\ld_0^{\m}} \Bigr)
   + \Bigl( \frac{c_0^2 \ld_4^{\m}}{ c_4^2 \ld_0^{\m}} \Bigr) 
     \Bigl(1 - \frac{2 \m}{\ld_4^{\m}} \Bigr) 
            = 12 \frac{1}{\sqrt{\m}} + O ( e^{-\sqrt{\m}/2} ),
\]
and we just get by.  Running through the same computation with the addition of
a test function of $p_{0,4}$ in the integrand completes the proof.

\section{Proof of the main error estimate}
\setcounter{num}{5}
\setcounter{equation}{0}
\label{sec:error}

Having renormalized in terms of the Gaussian measure $P_{\m}^*$ defined 
in Lemma {\ref{lem:starmeas}}, we now
return to the asymptotics of $f(\m, \ep)$, recall (\ref{equ:shift}).  
Using the form of $A(\cdot)$ we first write
\beq
\label{equ:reshuffle}
\lefteqn{   \sqrt{\frac{\pi}{2}}  \frac{1}{Z_{\m}^*} e^{-I_{\m}(p_{\m})}   f(\m; \ep) } \\
& = & 
   \int_0^1 \int_0^1 e^{2\int_y^x p_{\m}} E_{\m}^{*} \left[  e^{2  \int_x^y p - 2 \int_0^1 p_{\m} p^3  
                             - \frac{1}{2} \int_0^1 p^4 }
      R(p+p_{\m}), \  ||p||_{\infty} \le \ep \sqrt{\mu}  \right] \;dy\,dx \nonumber
\eeq
in which $Z_{\m}^*$ is the $P_{\m}^*$ mass,
$ Z_{\mu}^* = {E}_0^0 [ \exp{ \{ -\frac{1}{2} \int_0^1 ( q_{\m} -2\mu) p^2  \}} ] $.               
In these terms, 
the main result of this section is that right hand side of (\ref{equ:reshuffle}) equals $A(p_{\m}) R(p_{\m})$ up
to small multiplicative errors,
and, since $A(p_{\mu}) = \int_0^1 \int_0^1 e^{2 \int_y^x p_{\m}} dx dy$,
this is equivalent to the following.

\begin{thm}  For $\m \ra \infty$ and all $\ep > 0$ sufficiently small,
\[
   E_{\m}^{*} \Bigl [ \exp{  \Bigl\{ 2 \int_x^y p - 2 \int_0^1 p_{\m} p^3  
                             - \frac{1}{2} \int_0^1 p^4   \Bigr\} }
      R(p+p_{\m}), \  ||p||_{\infty} \le \ep \sqrt{\mu}  \Bigr] = R(p_{\m}) (1 + o(1)) 
\]
independently of $x, y \in [0,1]$.
\end{thm} 

This leaves the computation of $Z_{\m}^*$, taken up in the next section,
as the final ingredient in the proof of Theorem \ref{thm:main}.

As alluded to earlier, the fact that the (scaled) process $P_{\m}^*$ does not posses
a good spectral gap means we have to take extra care in dealing with the low modes.  In particular,
everywhere $p$ is split ($p = p_{\l} + p_{h}$) and
expanded, with differing considerations for each piece. To explain some of the difficulty, 
we introduce the shorthand 
\[
    p_{\l}(x) \equiv a_{\m}(x) \, \g_0 + b_{\m}(x) \, \g_2 +  c_{\m}(x) \, \g_3 
\]
and note that for large $\m$,  
$a_{\m}(x) \simeq (1/8) \cn^2(\sqrt{\m} x) +  O({\m}^{-1/2})$, 
$b_{\m}(x) \simeq  (1/2) \m^{-1/4} \sn(\sqrt{\m} x) \dn(\sqrt{\m} x)$,
and 
$c_{\m}(x) \simeq  (1/2) \m^{-1/4} \sn(\sqrt{\m} x) \cn(\sqrt{\m} x)$,
see again (\ref{equ:pexp}).  The point is that
while any positive moment of $|a_{\m}|$  decays like $\m^{-1/2}$, this mode remains $O(1)$
in sup-norm as $\m \ra \infty$.  The decay of $b_{\m}$ and $c_{\m}$ is only
slightly better, both of order $\mu^{1/4}$ in $L^{\infty}$ and order
$\mu^{-3/4}$ in $L^1$.  
It is not surprising then that we have to rely on 
cancellations in the lower modes.  For example, for the cubic we find that
\be
\label{equ:cancel}
 \int_0^1 p_{\m} p_{\l}^3 = 3 \Bigl( \int_0^1 p_{\m} a_{\m}^2 b_{\m} \Bigr)  \, \g_0^2 \g_2 
                    +  \int_0^1 p_{\m}( b_{\m} \g_2 + c_{\m} \g_3)^2, 
\ee
the potentially troublesome terms $\int_0^1 p_{\m} a_{\m}^3$ and $\int_0^1 p_{\m} a_{\m}^2 c_{\m}$
both vanishing. Even still, the remaining term $ \int_0^1 p_{\m} a_{\m}^2 b_{\m} 
\g_0^2 \, \g_2$ cannot be made small
on its own.   Besides cancellations, we need help from the (negative) quartic.  

Towards this end, we make the decomposition
\[ 
     2 \int_x^y p + 2  \int_0^1 p_{\m} p^3 - \frac{1}{2} \int_0^1 p^4  \equiv F_0(p,\m)
                        + F_1(p, \m)
\]
in which
\be
\label{equ:F0}
   F_0(p, \m)  =   6 \Bigl( \int_0^1 p_{\m} a_{\m}^2 b_{\m} \g_0^2 \g_2 \Bigr)  
                            - \frac{1}{2} \int_0^1 \Bigl( a_{\m}^4 \g_0^4 + b_{\m}^4 \g_2^4 + c_{\m}^4 \g_3^4 
                          \Bigr)
                       - \frac{1}{2} \int_0^1 p_h^4 - \int_0^1 p_{\l}^2 p_{h}^2  
\ee
and 
\beq
\label{equ:F1}
    F_1(p, \m) & = & 2 \int_0^1 p_{\m} p_h^3  + 6 \int_0^1 p_{\m} p_{\l} p_h^2 + 6 \int_0^1 p_{\m} p_{\l}^2 p_h \\
               &  & + 2 \int_x^y p  - 3 \int_0^1 p_{\l}^3 p_h  - 3 \int_0^1 p_{\l} p_{h}^3  \nonumber \\
               &  & + 2 \int_0^1 p_{\m} ( b_{\m} \g_2 + c_{\m} \g_3)^2 
                        - \frac{1}{2} 
                         \int_0^1 ( p_{\l}^4 - a_{\m}^4 \g_0^4 - b_{\m}^4 \g_2^4 - c_{\m}^4 \g_3^4 ). \nonumber
\eeq
With this in hand, we note that for any 
event ${\cal A}$ of the path:
\beqn
\lefteqn{
  \Bigl|  \, E_{\m}^* \Bigl[  e^{ F_0(p,\m) + F_1(p, \m) }  \,                     
      R(p+p_{\m}), \, {\cal A}  \Bigr] -  R(p_{\m})  E_{\m}^* \Bigl[ e^{F_0(p, \m)}, \, {\cal A} \Bigr] \, \Bigr| } \\
& \le  &  R(p_{\m})  
      E_{\m}^* \Bigl[  e^{F_0(p, \m)}  | F_1(p, \m) | e^{ |   F_1(p, \m) |} , \, 
            {\cal A}  \Bigr]
        +    E_{\m}^* \Bigl[ |  R(p_{\m}) - R(p+p_{\m}) |  e^{F_0(p,\m)} ,  
             \, {\cal A}    \Bigr]  \\
&   &   +   E_{\m}^* \Bigl[ | R(p_{\m}) - R(p+p_{\m})  | e^{F_0(p, \m)} \, | F_1(p, \m) | 
                         e^{ | F_1(p, \m) | }, 
 \, {\cal A}  \Bigr] 
\eeqn
where the elementary inequality  $|1 - e^f|$ $\le |f| e^{|f|}$ 
has been used.  Thus, after successive applications of H\"older's inequality, the
theorem is a consequence of the following three facts.

\begin{lemma}  
\label{lem:explicit}
There exists a $\theta^{\prime} > 1$ such that
\[
        E_{\m}^* \Bigl[ \exp \Bigl\{ \theta \, F_0(p, \m) \Bigr\} , \, ||p||_{\infty} \le \ep \sqrt{\m} \Bigr]
                = 1 + O(\m^{-1/2})
\]
for all $\theta \in [1, \theta^{\prime})$ as $\m \ra \infty$.
\end{lemma}

\begin{lemma} 
\label{lem:pain}
Given $\theta > 1$,
\[
  \lim_{\m \ra \infty}   E_{\m}^* \Bigl[   | F_1(p, \m)|^{\theta} \, 
 \exp{ \Bigl\{ \theta |F_1(p,\m)|} \Bigr\}  , \, 
|| p||_{\infty} \le \ep \sqrt{\m} \Bigr]  = 0.
\]
for all $\ep > 0$ small enough.
\end{lemma}

\begin{lemma}
\label{lem:R}
For any $\theta > 1$ there is the bound 
\[
     \Bigl( E_{\m}^* 
        \Bigl[ \Bigl|  R(p_{\m} + p ) - R(p_{\m}) \Bigr|^{\theta}, 
            \, ||p||_{\infty} \le \ep \sqrt{\m} \Bigr]  \Bigr)^{1/\theta} 
    \le C_0 \, R(p_{\m}) \, {\m}^{-3/4}
\]
in which $C_0$ depends on $\ep$ but not $\m$.
\end{lemma}

Furthermore, the proofs of these three ingredients 
rely to varying degree on the next lemma which describes
the convergence of $P_{\m}^*$ to the zero path as $\m \ra \infty$.

\begin{lemma}  
\label{lem:meansqdec} 
Under $P_{\m}^*$ the path has the following decay in mean-square.  We have 
\[
E_{\m}^* \Bigl[ \int_0^1 p^2(x) dx \Bigr] \le C_1 \, \m^{-1/2},
\]  
and, if we remove the low modes:
\[
\sup_{0 \le x \le 1} E_{\m}^* \Bigl[ p_h^2(x) \Bigr] \le C_2 \, \m^{-1/2}.
\]
\end{lemma}

\noindent{{\bf  Proof }  For the integrated 
mean-square estimate over the full path $p$ one just computes:
\be
\label{equ:vars}
   E_{\m}^* \Bigl[ \int_0^1 p^2(x) dx \Bigr]  =  \int_0^1 a_{\m}^2(x)dx + \frac{1}{\ld_{2}^{\m}- 2\m} 
                                                 + \frac{1}{\ld_{3}^{\m}- 2\m}  
                                  + \sum_{\l = 5}^{\infty} \frac{1}{\ld_{\l}^{\mu} - 2 \m}. 
\ee
The first three terms we know explicitly.  The $\m^{-1/2}$
decay of $\int_0^1 a_{\m}^2$ has already been remarked upon,
and  both $\ld_2^{\m}$ and  $\ld_3^{\m}$ are approximately  $5 \m$
for $\mu$ large.  We also know that
 $\ld_{\l}^{\m} - 2 \m = O(\m)$ for any fixed $\l \ge 2$ 
and $\m \ra \infty$. However, controlling the whole sum requires
more precise eigenvalue asymptotics provided by
the following classical result:
\be
\label{equ:eigexpand}
   \ld_{\l}^{\m} = 4 \pi^2 \l^2 + \int_0^1 {q_{\m}}  
+ O \Bigl( \frac{1}{4 \pi^2 \l^2} \int_0^1 (q_{\m} - {\int_0^1 q_{\m}})^2 \Bigr)
                 = 4 \pi^2 \l^2 + 6 \m + O \Bigl( \m \l^{-2} \Bigr)
\ee
for all large values of the index.  See, for example, Theorem 2.12 of \cite{MW}.
The tail of the sum in (\ref{equ:vars}) then behaves like
\[
   \sum_{\l \ge L} \frac{1}{\m + \l^2} \simeq  
\frac{1}{\sqrt{\m}} \int_{L/\sqrt{\m}}^{\infty} \frac{dx}{1 + x^2} = O(\m^{-1/2}),
\]
completing the verification.

As for
\[
 \sup_{0 \le x \le 1} E_{\m}^*[ p_h^2(x)] = \sup_{0 \le x \le 1} \, \Bigl\{ 
\sum_{\l = 5}^{\infty} \frac{(\phi_{\l}^{\m}(x))^2}{\ld_{\l}^{\m} - 2 \m}   \Bigr\},
\]
the result will follow granted a uniform bound 
on $||\phi_{\l}^{\m}||_{\infty}$ independent of $\m$ and of the index $\l \ge 5$. 
(The intuition is thus: $\l \ge 4$ corresponds
to continuous spectrum for $\m \ra \infty$, and so the corresponding
eigenfunctions should 
remain ``flat'' for large values of $\m$.) 
Once more, the ingredients of the verification are classical.  It is convenient to consider 
the unscaled equation which reads $u^{\pr \pr}(x) + ( 6 k^2 \cn^2(x) + \gamma_{\l}) u(x) = 0 $
($u(0) = u(\sqrt{\m})$) with $\gamma_{\l}> 0$ and of order
$\l^2/\sqrt{\m}$ (note again \ref{equ:eigexpand}).  Then, for any $\l$ up
to order $\sqrt{\m}$, a comparison argument with the explicit  $\l = 4$ case ($\gamma \simeq 0$) 
will produce the desired bound.  On the other hand, when $\l > C \sqrt{\m}$ with $C$ large, well
known arguments (see \cite{Hoch1} for a model) will show that the solution
is uniformly approximated (up to errors of $O(\l^{-1})$) by a single trigonometric function with $\l$ 
oscillations. The proof is complete.

\bigskip

\noindent{{\bf Proof of Lemma {\ref{lem:explicit}}} }  We concern ourselves with the upper bound, the lower
bound being a simple consequence of Jensen's inequality.  
First note that
\[
   E_{\m}^{*} \Bigl[ e^{ \theta \, F_0(p, \m)}, \, ||p||_{\infty} \le \ep \sqrt{\m} \Bigr]
   \le  E \Bigl[ \exp \Bigl\{ 6 \, \theta \, \Bigl( \int_0^1 p_{\m} a_{\m}^2 b_{\m} \Bigr)
             \,  \g_0^2 \g_2  - \frac{\theta}{2} 
                           \Bigl( \int_0^1
                              a_{\m}^4 \Bigr) \, \g_0^4 \Bigr\} \Bigr]. 
\] 
That is, we only really use the quartic in the lowest mode to control the bad cubic term. 
Next, the integral on the right hand side may be performed in the $\g_2$
variable, and we 
will be satisfied to show that there exists a choice of $\theta > 1$ and a  $\d > 0$ such that
\be
\label{equ:balance} 
    36 \; \theta^2 \, \Bigl( \int_0^1 p_{\m} a_{\m}^2(x) b_{\m}(x)  \, dx \Bigr)^2 
                   -  \theta \, \int_0^1 a_{\m}^4(x) dx   \le  - \frac{1}{\sqrt{\m}} \d 
\ee
for all large enough $\m$. Now 
\[
    \Bigl( \int_0^1 p_{\m}(x) a_{\m}^2(x) b_{\m}(x) dx \Bigr)^2 
                      = \frac{1}{\sqrt{\m}}   \frac{1}{8^4} 
                       \frac{ ( \int_0^{\sqrt{\m}} \cn^4(x) \sn^2(x) \dn(x) dx )^2 }
                            { 3  \int_0^{\sqrt{\m}} \sn^2(x) \dn^2(x) dx }                           
\] 
and
\[
    \int_0^1  a_{\m}^4(x) dx =   \frac{1}{\sqrt{\m}}   \frac{1}{8^4}  \int_0^{\sqrt{\m}} \cn^8(x) dx  
\]
up to (unimportant) errors of order $\m^{-1}$.  
Further, up to errors exponentially small in $\sqrt{\m}$,
the integrals on the right of the last 
two displays may be replaced by integrals over the whole line of the corresponding
hyperbolic-trigonometric functions.  That is, the validity of the desired
inequality
(\ref{equ:balance}) is equivalent to whether
\[
    12 \int_{-\infty}^{\infty} \sech^{5}(x) \tanh^2(x) dx 
   < \int_{-\infty}^{\infty}  \sech^{8}(x) dx  \, \int_{-\infty}^{\infty} 
                  \sech^2(x) \tanh^2(x) dx. 
\]
Both sides of the latter may be worked out to read $12 (\pi/16)^2 < (2/3)
 (32/35)$ 
which is indeed true.  The proof is finished.

\bigskip

\noindent{{\bf Proof of Lemma {\ref{lem:pain}}} }  Consider the first line
comprising $F_1(p, \m)$, see (\ref{equ:F1}).  
That is, the estimate is detailed for $\tilde{F}_1(p,\m) 
\equiv 2 | \int_0^1 p_{\m} p_h^3 | + 6 | \int_0^1 p_{\m} p_{\l} p_h^2 |
+ 6 |  \int_0^1 p_{\m} p_{\l}^2  p_h | $. These terms are in a sense 
the most difficult as they involve the additional factor of $\sqrt{\m}$
through $p_{\m}$.  At the end we comment on how to deal with the remaining
terms in the full $F_1$.

For the present task, it suffices to prove the following two types of
estimate.  First,
\be 
\label{equ:means}
 E_{\m}^* \Bigl[  \Bigl| \int_0^1 p_{\m} p_{h}^3  \Bigr|^m \Bigr] \, \ra \, 0, \; \;
 E_{\m}^* \Bigl[  \Bigl| \int_0^1 p_{\m} p_{\l} p_{h}^2  \Bigr|^m \Bigr] \, \ra \, 0, \; \;  
  {\mbox{ and }} \;  E_{\m}^* \Bigl[  \Bigl|  \int_0^1 p_{\m} p_{\l}^2  p_h  \Bigr|^m  \Bigr] \, \ra \, 0 .
\ee
as $\m \ra \infty$ for all $m$ large enough. 
Second, for whatever $C > 0$
\beq
\label{equ:expbound} 
\lefteqn{   \limsup_{\m \ra \infty}    
 \, E_{\m}^* \Bigl[ \exp{ \Bigl\{ C \sqrt{\m} \int_0^1 |p_h|^3  \Bigr\} }, ||p||_{\infty} \le \ep \sqrt{\m} \Bigr] 
         < \infty,} \\
&  &             
   \limsup_{\m \ra \infty}  \, E_{\m}^* \Bigl[ \exp{ \Bigl\{ C \sqrt{\m}  \int_0^1 |p_{\l}| p_h^2  \Bigr\}},
                            \, ||p||_{\infty} \le \ep \sqrt{\m} \Bigr] < \infty,  \nonumber \\ 
&  & {\mbox{ and } }
\limsup_{\m \ra \infty} \,  E_{\m}^* \Bigl[ \exp{ \Bigl\{ C | \int_0^1 p_{\m} p_{\l}^2  p_h |  \Bigr\}},
                            \, ||p||_{\infty} \le \ep \sqrt{\m} \Bigr] < \infty \nonumber 
\eeq
when $\ep > 0$ is chosen appropriately.

Starting with (\ref{equ:means}) and working left to right we first have
\beqn
   E_{\m}^* \Bigl[  \Bigl| \int_0^1 p_{\m}(x) p_h^3(x) dx  \Bigr|^{m} \Bigr]
   & \le & \mu^{m/2}  \int_0^1 E_{\m}^*[ |p_h|^{3m}(x) ] dx \\
   &  =  & C_3 \,  \mu^{m/2} \int_0^1 \Bigl(  E_{\m}^*[ p_h^{2}(x) ] \Bigr)^{3m/2}  dx \le C_4 \, \m^{-m/4}.
\eeqn
Here we have used Jensen's inequality, the fact that $p_h(x)$ is Gaussian, and Lemma 
\ref{lem:meansqdec}.  In a similar fashion
\beqn 
\lefteqn{ E_{\m}^* \Bigl[  \Bigl( \int_0^1 |p_{\m}|(x) |p_{\l}| p_h^3(x) dx  \Bigr)^{m} \Bigr]} \\
& \le & C_4 \,\m^{m/2} \,   
  \Bigl\{ E_{\m}^* \Bigl[  \Bigl( \int_0^1 |a_{\m}|(x) p_h^2(x) dx  \Bigl)^m 
         +       \Bigl( \int_0^1 |b_{\m}|(x) p_h^2(x) dx  \Bigl)^m  +  
               \Bigl( \int_0^1 |c_{\m}|(x) p_h^2(x) dx  \Bigl)^m  \Bigr\}.
\eeqn
Restricting attention to the first term as $a_{\m}$ has less decay compared with 
$b_{\m}$ or $c_{\m}$ we find that by Jensen's inequality and Lemma \ref{lem:meansqdec},
\[
\m^{m/2}  \, E_{\m}^* \Bigl[  \Bigl( \int_0^1 |a_{\m}|(x) p_h^2(x) dx  \Bigl)^m  \le  C_5
   \m^{m/2}  \Bigl( \int_0^1 |a_{\m}|^{m}(x) dx \Bigr) \Bigl( \sup_{0 \le x \le 1}   
              E_{\m}^* [ p_h(x) ]  \Bigr)^{m/2} \le C_6 \, \m^{-1/2}.
\]
As for the last expectation in (\ref{equ:means}), it may be bounded by  
\[
\m^{m/2}  E_{\m}^* \Bigl[ \Bigl( \int_0^1    p_{\l}^2(x)  |p_h (x) | \, dx \Bigr)^m  \Bigr] 
     \le   C_7 \, \m^{m/2}  \, \Bigl\{   E_{\m}^* \Bigl[ \Bigl( \int_0^1  a_{\m}^2(x) |p_h (x) | \, dx 
                   \Bigr)^m  \Bigr] + {\mbox{ like terms in }} b_{\m}^2 {\mbox{ and }}  c_{\m}^2 \Bigr\} .
\] 
Spelling out the first term involving we find 
\beqn
\lefteqn{  \m^{m/2} \, E_{\m}^* \Bigl[ \Bigl( \int_0^1  a_{\m}^2(x) |p_h (x) | \, dx \Bigr)^m  \Bigr]} \\
   & \le &  C_8 \, \m^{m/2} \,
\Bigl( \int_0^1 |a_{\m}|^{\frac{2 m}{m-1}}(x) dx  \Bigr)^{m-1}   \, 
              \int_0^1 \Bigl( E_{\m}^* [  p_h^2(x) ] \Bigr)^{m/2} dx 
\le C_9 \m^{- \frac{(m-1)}{2}}
\eeqn
after an application of   H\"older's inequality, and, again, the $\m^{-1/2}$ decay of 
the moments of $a_{\m}$ along with Lemma \ref{lem:meansqdec}.  Terms two
and three are handled in the same way.

Turning now to exponential bounds ({\ref{equ:expbound}), we  require
two additional observations.  
The first is that there exists an $M \ge 1$ so that 
\be
\label{equ:splitbound}
   ||p||_{\infty} \le \ep \sqrt{\m}   \; \; \; {\mbox{ implies }}  \; \; 
        || p_{\l} ||_{\infty} 
            \le  M \ep \sqrt{\m};  \; {\mbox{ and }} \; \; 
    ||p_h||_{\infty} \le  M  \ep \sqrt{\m}. 
\ee
for some $M > 0$.  To see this, first multiply the  inequalities 
$-\ep \sqrt{\m} \le p(x) \le - \ep \sqrt{\m} $ through by $\phi_0^{\m}(x) > 0$
to find that
\[
    | \g_0 | \,  \int_0^1 a_{\m}(x)  \phi_0^{\m}(x) dx  \le \ep \sqrt{\mu} \, \int_0^1 \phi_0^{\m}(x) dx.  
\]
Next it may be checked that $\int_0^1 a_{\m}  \phi_0^{\m} $ 
and $\int_0^1 \phi_0^{\m}$ are of the same order, and so 
\be
\label{equ:g0}
                 | \g_0 | \le C_{10} \, \ep \m^{1/2}.
\ee
The verification of (\ref{equ:splitbound}) is completed by showing that
\be
\label{equ:g23}
        |\g_2| \le C_{11} \, \ep \m^{3/4}   \; \; {\mbox{ and }}  \; \;  |\g_3| \le  C_{12} \,  \ep \m^{3/4} 
\ee
on $||p||_{\infty} \le \ep \sqrt{\mu}$.  This is similar: the inequality 
is now multiplied through by $( ||\phi_2^{\m}||_{\infty} + \phi_2^{\m})$ 
( or   
$(||\phi_3^{\m}||_{\infty} + \phi_3^{\m})$) and integrated. Since the
path is mean zero this produces an inequality in $|\g_2|$ (or $|\g_3|$) alone which 
is equivalent to (\ref{equ:g23}).  The second fact we will need is the following upper
bound on the sup-norm deviations of $p_h$.  For all $N$ large, there exist (positive) 
constants $a$ and $b$ independent of $\m$ such that
\be
\label{equ:borell}
 P_{\m}^* \Bigl( || p_h ||_{\infty} > N \Bigr) \le  a \, \exp \Bigl[ - b \sqrt{\m} N^2  \Bigr].   
\ee 
This is an immediate consequence of Borell's inequality (see \cite{A}) and Lemma \ref{lem:meansqdec}.

With those points established, we turn to the first expectation in question:
\beq
\label{equ:firstexp}
\lefteqn{   E_{\m}^* \Bigl[ \exp{ \Bigl\{ C \sqrt{\mu} \int_0^1 |p_h|^3 \Bigr\}}, \, 
                        ||p||_{\infty} \le \ep \sqrt{\m} \Bigr] } \\
& \le &    E_{\m}^* \Bigl[ \exp{ \Bigl\{  C \sqrt{\mu} N \int_0^1 p_h^2  \Bigr\}}, \, ||p_h||_{\infty} \le N \Bigr] +
            E_{\m}^* \Bigl[ \exp{ \Bigl\{  C \mu \ep \int_0^1 p_h^2  \Bigr\} }, \, 
||p_h||_{\infty} \ge N \Bigr] \nonumber \\
& \le &  E_{\m}^* \Bigl[ \exp{ \Bigl\{  C \sqrt{\mu} N \int_0^1 p_h^2  \Bigr\}} \Bigr] +
             \Bigl( E_{\m}^* \Bigl[ \exp{ \Bigl\{  2 C \mu \ep \int_0^1 p_h^2 \Bigr\}} \Bigr] \Bigl)^{1/2}  
            \Bigr(  P (  ||p_h||_{\infty} \ge N ) \Bigr)^{1/2}  
                   \nonumber \\
& \le &  \exp \Bigl[ C N \cdot \sum_{\l \ge  5 } \frac{\sqrt{\m}}{ \ld_{\l}^{\m} - 2 \m} \Bigr]
          +    a \exp \Bigl[ C  \ep \cdot  \sum_{\l \ge  5 } \frac{\m}{ \ld_{\l}^{\m} - 2 \m} \Bigr] 
             \, \exp \Bigl[  - \frac{1}{2} b \sqrt{\m} N^2 \Bigr]. \nonumber
\eeq
Note that this last inequality holds as soon as $\m$ is large enough and $\ep$ small enough so that
both $ C N \cdot (\sqrt{\m} / (\ld_{2}^{\m} - 2 \m)) \simeq (1/3) C N \mu^{-1/2} $ and 
$ 2 C \ep \cdot ( \sqrt{\m}/ (\ld_{2}^{\m} - 2 \m))  \simeq (2/3) C \ep $ are less than one-half.  Now we recall 
that the sum $ \sqrt{\m} \sum_{\l > 2} \ld_{\l}^{\m}$ is bounded by a fixed constant for $\m \up \infty$.
Indeed, this is really the content of Lemma \ref{lem:meansqdec},  this sum being finite along with
$\limsup_{\m \ra \infty} \sum \frac{\sqrt{\m}}{\m + k^2} < \infty$.  So the first term in the last line above is bounded independent of $N$,
and the second will go to zero as $\m \ra \infty$ by appropriate choice of $N$.

The next integral is approached the same way.  We find that
\beqn
\lefteqn{ E_{\m}^* \Bigl[ \exp{ \Bigl\{ C \sqrt{\mu}  
\int_0^1  |p_{\l}| p_h^2 \Bigr\}}, \, ||p||_{\infty} \le \ep \sqrt{\m} \Bigr] }  \\
& \le &    E_{\m}^* \Bigl[ \exp{ \Bigl\{ C\sqrt{\m} N^2 \Bigl( ||a_{\m}||_1  |\g_0|  + ||b_{\m}||_1 | \g_2 |
                                        +  || c_{\m} ||_1 | \g_3 | \Bigr)    \Bigr\}} \Bigr]  \\
&     &    +
        E_{\m}^* \Bigl[ \exp{ \Bigl\{  C \cdot M \ep \mu  \int_0^1 p_h^2 \Bigr\}}, \,     
     ||p_h||_{\infty} \ge N \Bigr]. 
\eeqn
In the first term on the right hand side we have used  $||p_h||_{\infty} \le N$ by explicitly restricting to that
set of paths.   Since  $ || a_{\m}||_1 = O(\mu^{-1/2})$ and $|| b_{\m}||_1, || c_{\m}||_1 = O(\m^{-3/4})$,  
this term is plainly bounded.
In the second term we have used (\ref{equ:splitbound}): 
$|p_{\l}| \le M \ep \sqrt{\mu}$ due the overall control of $||p||_{\infty}$. 
From this point this term yields to considerations identical 
to those for the second term in line of (\ref{equ:firstexp}).

Last we come to bound the expectation involving a large constant multiple of
\beq
\label{equ:anotherexp}
\lefteqn{  \Bigl| \int_0^1 p_{\m}(x) p_{\l}^2(x) p_h(x) dx \Bigr|
    \le    \Bigl| \int_0^1 p_{\m}(x) a_{\m}^2(x) p_h(x) dx \Bigr| \, \g_0^2}   \\
   &     & +  \,     \Bigl| \int_0^1 p_{\m}(x) b_{\m}^2(x) p_h(x) dx \Bigr| \, \g_2^2 + 
         2  \, \Bigl| \int_0^1 p_{\m}(x) a_{\m}(x) b_{\m}(x) p_h(x) dx \Bigr| \, |\g_0 \g_2| + etc; \nonumber
\eeq
the $etc$ indicating like terms in $\g_3^2$, $|\g_0 \g_3|$ and $| \g_2 \g_3 |$.  Making use of (\ref{equ:g0})
and (\ref{equ:g23}) we note that, for example, 
\[
\Bigl| \int_0^1 p_{\m}(x) a_{\m}(x) b_{\m}(x) p_h(x) dx \Bigr|  | \g_0 \g_2 |
 \le C_{13} \, \ep \m  \,
 \Bigl| \int_0^1 \sn(\sqrt{\m} x) a_{\m}(x) {\tilde \phi}_2^{\m}(x) p_h(x) dx \Bigr|   \, | \g_0 | .
\]
Similarly, each term on the right
hand side of (\ref{equ:anotherexp}) may  be bounded in turn
by (constant multiples of) expressions of the form 
\[
\m \ep \Bigl| \int_0^1 \psi(\sqrt{\m} x) p_h(x) dx \Bigr| |\g_{\l}| 
      \equiv \m \ep | {\cal G}({\psi}) | |\g_{\l} |
\]
with $\l = 0,2,$ or $3$ and  $\psi$ a smooth periodic function 
for which $\limsup_{\m \ra \infty} \int_0^{\sqrt{\m}} |\psi(x)| dx < \infty$.  
While $\psi$ differs in each appearance it now suffices to show that
for any such $\psi$, 
$\ep$ may be chosen small enough so that 
\[
  E  \Bigl[  \exp \Bigl\{  \m \ep | {\cal G}(\phi) | |\g_0| \Bigr\} \Bigr] 
  \le  C_{14}      E  \Bigl[  \exp \Bigl\{  \m^2 \ep^2 {\cal G}^2(\phi)  \Bigr\} \Bigr] 
 \le C_{15}  
\]  
independently of $\m$.  Now since ${\cal G}(\phi)$ is itself a mean zero Gaussian random
variable,  this last display is true as long as  
\be
\label{equ:Gvar}
      \m^2 E  [ {\cal G}^2(\psi)]  
= \m^2 \sum_{\l \ge 5} \frac{1}{\ld_{\l}^{\m} - 2 \m } 
 \Bigl( \int_0^1 \psi(\sqrt{\m} x) \phi_{\l}^{\m}(x) dx \Bigr)^2
\ee
may be bounded independently of $\m$. In this direction, we have the estimate
\be
\label{equ:RL}
     \int_0^1 \psi(\sqrt{\m} x) \phi_{\l}^{\m}(x) dx \le C_{15}  \Bigl(  \frac{1}{\sqrt{\m}} \wedge \frac{1}{\l}  \Bigr)
\ee
which is explained as follows.  In Lemma \ref{lem:meansqdec} it is remarked that there is a uniform $L^{\infty}$
bound on the high eigenfunctions $| \phi_{\l}^{\m}|$.  The $\m^{-1/2}$ decay may then always be extracted 
from the same decay of $\psi$ in $L^{1}$.  On the other hand, as also remarked in Lemma 
\ref{lem:meansqdec},
for high values of the index $\phi_{\l}^{\m}$ is well approximated by a trigonometric function of $\l$-turns.
The $\l^{-1}$ factor then follows from the lemma of Riemann-Lebesgue.
Finally, substituting (\ref{equ:RL}) into (\ref{equ:Gvar})  produces
\[
  \m^2 E   [ {\cal G}^2(\psi)]  \le C_{16} \Bigl(    \m^2 \sum_{4 \le \l  \le  \sqrt{\m}}  \frac{1}{\m + \l^2} \cdot \frac{1}{\m} + 
  \m^2 \sum_{\sqrt{\m} \le \l < \infty}  \frac{1}{\m + \l^2} \cdot \frac{1}{\l^2} \Bigr), 
\]
and it is readily checked that the latter remains bounded as  $\m \ra \infty$.

We close with some comments regarding the remaining terms in $F_1(p, \m)$. Looking at line 
two of (\ref{equ:F1}) it is plain at this point that the linear 
term $|\int_x^y p | \le \int_0^1 |p|$ poses no difficulty.  On the other hand 
$\int_0^1 |p_{\l}| |p_h|^3 \le M \ep \sqrt{\m} \int_0^1 |p_h|^3$ on the domain  
of integration which reduces the term to a type already dealt with.  Also, expanding the low
modes out in $| \int_0^1 p_{\l}^3 p_h |$ and employing (\ref{equ:g0}) and (\ref{equ:g23})
make this term comparable to $| \int_0^1 p_{\m} p_{\l}^2 p_h |$.  Finally, the last line 
in (\ref{equ:F1}) is explicit in terms of the low modes.  It may be done by hand,  cancellations
in the cross terms (as in (\ref{equ:cancel})) being of great importance.  The tedious details are not reported.
The proof is complete.

\bigskip

\noindent{{\bf Proof of Lemma {\ref{lem:R}}} } 
We actually prove that, on the set of paths $\{ p : ||p ||_{\infty} \le \ep \sqrt{\m} \}$,
there is a fixed constant so that
\[
   \Bigl| R(p_{\m} + p) - R(p_{\m}) \Bigr| \le C_0 \, R(p_{\m}) \, \frac{1}{\sqrt{\m}} \int_0^1 |p|.
\]
The $L^{\theta}$ estimate for $\theta > 1$ then follows from Jensen's
Inequality and the first conclusion of Lemma  \ref{lem:meansqdec}.

Introduce the un-normalized extremum $p^*(x) = k \, \sn(\sqrt{\m} x, k)$ 
and the functional
\[
   {\cal R}(a,p)  = \int_0^1 {\tilde \phi}_1^{\m}(x + a) p(x)  dx = \int_0^1 \cn \dn ( \sqrt{\m} (x + a) ) p(x) dx.
\]
It is the
derivative of ${\cal R}$ in $a$ which figures in the definition of $R$.
Clearly, for any path $p$, ${\cal R}(a,p)$  is 
is an analytic function of $a$, and  ${\cal R}^{\pr}(a, p) \le C
||p||_{\infty}$ with a constant independent of $a$.
Further, at the extremum $p^*$, ${\cal R}(a, p^*)$ has exactly two zeros at
$a = 0$ and $a = 1/2$, and the derivatives at those points satisfy  
$ {\cal R}^{\prime}(0, p^*) =$ 
$ -  {\cal R}^{\prime}(1/2, p^*)$. Also, one may check that $| {\cal
  R}^{\prime}(0, p^*)|$  
converges to a  positive constant as $\m \up \infty$.  
It follows that if $p$ satisfies $||p||_{\infty} \le \sqrt{\m} \ep$,  then,
for all large $\m$, ${\cal R}(a,p^*+p/\sqrt{\m})$ 
has exactly two zeros,  $r_1= r_1(p)$ and $r_2= r_2(p)$, lying within
neighborhoods 
of radius $\ep$ from $0$ and $1/2$ respectively.
In fact, $|r_1|$ and $|r_2 - 1/2|$ 
may be bounded in a more useful way by a constant multiple of $\frac{1}{\sqrt{\m}} \int_0^1 |p|$.  
Take the case of $r_1$:
\beq
\label{equ:r1b}
  \frac{1}{\sqrt{\m}} \int_0^1 |p(x)| dx  & \ge &    
\frac{1}{\sqrt{\m}}  \Bigl| \int_0^1 \tilde{\phi}_1^{\m}(x+ { r}_1) p(x) \, dx  \Bigr| \\
        &  = & \Bigl| \int_0^1 \Bigl( \tilde{\phi}_1^{\m}(x) - \tilde{\phi}_1^{\m}(x+ {r}_1)  \Bigr) p^*(x) \, dx  \Bigr| 
         =   | { r}_1(p)   | \,  \Bigl| {\cal R}^{\prime}(\tilde{r}, p^*)  \Bigr|  \nonumber
\eeq
with some $ \tilde{r} $ between $0$ and  $r_1$, that is, 
$|\tilde{r}| = O(\ep)$.  Since  there exists a $\delta > 0$ with $|  {\cal
  R}^{\prime}(0, p^*) | > \delta$ for all large $\m$, by analyticity 
$|  {\cal R}^{\prime}({\tilde r}, p^*) |$ satisfies a similar lower bound.  The claimed bound on $r_1$ follows.  

Next let us write ${\cal R}_0 = | {\cal R}^{\pr}(0, p^*) |$, ${\cal R}_1  = |{\cal R}^{\pr}(r_1, p^* + p/\sqrt{\m}) |$,
and ${\cal R}_2  = |{\cal R}^{\pr}(r_2, p^* + p/\sqrt{\m}) |$ in terms of which 
\[
\Bigl|  R(p_{\m} + p) - R(p_{\m}) \Bigr|  
=  \frac{\sqrt{\m}}{\sqrt{\int_0^1 (\tilde{\phi}_1^{\m})^2}} \;
         \Bigl| \frac{{\cal R}_0}{2}   -  \frac{ {\cal R}_1 {\cal R}_2}{
	   {\cal R}_1 + {\cal R}_2} \Bigr|
\le  \frac{\sqrt{\m}}{\sqrt{\int_0^1 (\tilde{\phi}_1^{\m})^2}} \;
        \left\{ \Bigl|  \frac{ {\cal R}_1 ({\cal R}_0- {\cal R}_2)}{ {\cal
	    R}_1 + {\cal R}_2} \Bigr| 
+\Bigl| \frac{ {\cal R}_2 ({\cal R}_0- {\cal R}_1)}{ {\cal R}_1 + {\cal R}_2} \Bigr| \right\}.    
\]
It follows that it enough to show that $| {\cal R}_0 - {\cal R}_1 |$ and $|{\cal R}_0 - {\cal R}_1 |$
are bounded by  $C_1 {\cal R}_0  \m^{-1/2} \int_0^1 |p|$ or, what is the same, by  $C_2 \int_0^1 |p| $.  
(That $ {\cal R}_{1,2} / (  {\cal R}_1 + {\cal R}_2 )  =  O(1)$ is plain.)  Consider the difference
of ${\cal R}_0$  and $ {\cal R}_1$, the other estimate being identical.  On the set 
$\{ p : ||p||_{\infty} \le \sqrt{\m} \, \ep \}$ we have
\beqn
\lefteqn{   | {\cal R}_0 - {\cal R}_1 |   } \\
& = & 
           \Bigl| \int_0^1   ( {\tilde \phi}_1^{\m})^{\pr} (x )   p^*(x) dx  
     - \int_0^1  ( {\tilde \phi}_1^{\m})^{\prime} (x+ r_1 ) \Bigl(p^*(x) + \frac{1}{\sqrt{\m}} p(x) \Bigr)  dx \Bigr|  \\
& \le &  \int_0^1   \Bigl|   ( {\tilde \phi}_1^{\m} )^{\pr} (x) -
                                    (  {\tilde \phi}_1^{\m} )^{\pr} (x + r_1) \Bigr|  | p^*(x) |   \, dx  + 
                                                 \frac{1}{\sqrt{\m}} 
                               \int_0^1  \Bigl| ( {\tilde \phi }_1^{\m} )^{\pr}(x + r_1) \Bigr|  |p(x)| dx \\
& \le &   \m \, | r_1(p) | \, \int_0^1 \Bigr|  (  \cn \dn )^{\pr \pr} ( \sqrt{\m} (x + r^*) ) \, \sn(\sqrt{\m} x) \Bigr|  \, dx
                                              +       \int_0^1 | (\cn \dn)^{\pr} ( \sqrt{m}(x + r_1)) | 
                                                            \,  |p(x)|  \, dx                                  
\eeqn 
for some $r^*$, $-|r_1| \le r^* \le |r_1|$.  In line two we have supposed that $\ep$ is small enough that 
the absolute values in the definition of ${\cal R}_0$ and ${\cal R}_1$ may be left off.  To finish, note 
that $ | (\cn \dn)^{\pr}(\cdot)| \le 2 $, 
\[
\int_0^1 \Bigr|  (  \cn \dn )^{\pr \pr} ( \sqrt{\m} (x + r^*) ) \, \sn(\sqrt{\m} x) \Bigr|  \, dx = O ( \m^{-1/2}) 
\]
if $ r^*$ is bounded with probability one, and last, from (\ref{equ:r1b}), $|r_1(p)| \le C_3
\m^{-1/2} \int_0^1 |p|$.  The proof is complete.

\section{The Gaussian correction}
\setcounter{num}{6}
\setcounter{equation}{0}
\label{sec:gaussian}

All would be for naught if we could not compute 
the ``Gaussian correction'':
\begin{equation}
\label{equ:Z}
 Z(\mu) = E_0^0  \Bigl[ 
            \exp\{ - \frac{1}{2} 
         \int_0^1 ( q_{\mu}(x)  - 2 \mu ) p^2 (x) dx          
                \} \Bigr] P_0 \Bigl(  \int_0^1 \phi_1^{\m}(x) p(x) dx = 0 \Bigr)  
\end{equation}
$( = Z^* P_0 (\int_0^1  \phi_1^{\m} p = 0))$.
The properties of the operator $Q_{\m} = - d^2/dx^2 +q_{\mu}$ outlined in Section \ref{sec:hillspec}
now play an essential role.  As in the proof of Lemma \ref{lem:starmeas}, 
it is convenient to renormalize and split
the $E_0^0$ integral into two pieces as follows
\beq
\label{equ:Zsplit}   
E_0^0 \Bigl[ \exp{ \Bigl\{ - \frac{1}{2} \int_0^1 (q_{\m} - 2 \m) p^2 \Bigr\} } \Bigr] & = &
 E_{q_{\m}} \Bigl[ \exp{ \Bigl\{ \mu \int_0^1 p^2  \Bigr\} } \  \Bigl| 
                        \ \int_0^1 p = 0, \int_0^1 \phi_1^{\m}  p = 0 \Bigr]  \\
&   &  \times  E_0^0 \Bigl[ \exp{ \Bigl\{ - \frac{1}{2} \int_0^1 q_{\mu} p^2  \Bigr\} } \Bigr]. 
                \nonumber
\eeq
Again, $E_{q_{\m}}$ denotes the mean-value with respect to the 
Gaussian weight 
$Z_{q_{\m}}^{-1} \exp{[ - 1/2 \int_0^1 q_{\mu} p^2 ]} \times d CBM$.
We recall that
the point of this splitting 
is that, without any conditioning, 
the latter is a proper probability measure on periodic paths ({\em i.e.},
$Z_{q_{\m}} < \infty$).   

Taking advantage of
this, the first integral in  (\ref{equ:Zsplit}) may be evaluated as follows.

\begin{prop} 
\label{lem:g1}
For all $\m > 0$, 
\beqn
\lefteqn{ 
   E_{q_{\m}} \Bigl[ \exp{ \Bigl\{ \mu \int_0^1 p^2 \Bigr\} 
    } \ \Bigl| \ \int_0^1 p  = 0, \int_0^1 \phi_1^{\m} p = 0 \Bigr] } \\
& = &  \Bigl[  \Bigl(  1 - \frac{2 \m }{\ld_0^{\m}} \Bigr)  
               \Bigl(  1 - \frac{2 \m }{\ld_1^{\m}} \Bigr) 
               \Bigl(  1 - \frac{2 \m }{\ld_4^{\m}} \Bigr) \Bigr]^{1/2} \times
\Bigl[  \frac{\ld_4^{\m} c_0^2 + \ld_0^{\m} c_4^2}{ (\ld_4^{\m} - 2\m) c_0^2 
                           + (\ld_0^{\m} - 2 \m) c_4^2 } \Bigr]^{1/2} 
\times \sqrt{\frac{\Delta^2(0) - 4}{\Delta^2(2 \mu) - 4}}.
\eeqn
Here $\Delta(\ld)$  is the discriminant of $Q_{\m}$ and the constants 
$c_0$ and $c_4$ are as defined in (\ref{equ:c's}).
\end{prop}

The computation behind the second piece is more involved.  The result is:

\begin{prop} 
\label{lem:g2}
There is the explicit formula 
\be
\label{equ:g2}
 E_0^0 \Bigl[ \exp{ \Bigl\{ - \frac{1}{2} \int_0^1 q_{\mu} p^2 
                    \Bigr\}   } \Bigr] P_0 \Bigl( \int_0^1 \phi_1^{\m} p  = 0  \Bigr)
  =   \Bigl[  \frac{2 \pi}{\ld_1^{\m}} ( \frac{c_0^2}{\ld_0^{\m}} 
                     + \frac{c_4^2}{\ld_4^{\m}} ) \Bigr]^{-1/2}  
               \times \frac{1}{\sqrt{ \Delta^2(0) - 4 } };   
\ee
the notation being the same as in the previous result.
\end{prop}

Note what has been accomplished:  by Propositions \ref{lem:g1} and {\ref{lem:g2}} 
and Hochstadt's formula (\ref{equ:hoch}), 
$Z(\m)$ is now expressed completely in terms of the simple spectrum of $Q_{\m}$
which we know explicitly.

\bigskip

\noindent{\bf Proof of Proposition \ref{lem:g1} }
Once more we bring in the expansion of the path
$
  p(x) = \sum_{\l=0}^{\infty} \frac{1}{\sqrt{\ld_{\l}^{\m}}} \phi_{\l}^{\m}(x) {\g}_{\l}
$
under $P_{q_{\m}}$.  This translates the expectation of interest to:
with $E$ the mean corresponding to the independent Gaussian $\g$'s, 
\beq
\label{equ:expand1}
\lefteqn{
E_{q_{\m}} \Bigl[ \exp{ \Bigl\{ \m \int_0^1 p^2(x) dx \Bigr\}  } \ \Bigl| 
 \ \int_0^1 p(x) dx = 0, \int_0^1 \phi_1^{\m}(x) p(x) dx = 0 \Bigr] } \\
& = & E \Bigl[   \exp{ \Bigl\{ \sum_{\l = 0}^{\infty} \frac{ \m }{\ld_{\l}^{\m}}  \g_{\l}^2  \Bigr\} }
                \  \Bigl| \ \g_1 = 0, \frac{c_0}{\sqrt{\ld_0^{\m}}}  
 \g_0 - \frac{c_4}{\sqrt{\ld_4^{\m}}} \g_4 = 0   \Bigr] \nonumber \\
& = & \Bigl[  \prod_{ 2 \le \l \le \infty, \l \neq 4} 
\Bigl( 1 - \frac{2 \m}{\ld_{\l}^{\m}} \Bigr)  \Bigr]^{-1/2}   
              E \Bigl[  \exp{ \Bigl\{ \frac{\m}{\ld_0^{\m}} \g_0^2 +   
 \frac{\m}{\ld_0^{\m}}  \g_4^2 \Bigr\} }  \Bigl| 
                           \frac{c_0}{\sqrt{\ld_0^{\m}}} \g_0 - 
 \frac{c_4}{\sqrt{\ld_4^{\m}}} \g_4 = 0   \Bigr]. \nonumber 
\eeq
The integral in the last line already appeared in the proof of Lemma \ref{lem:starmeas};  
we now detail its evaluation:
\beqn
\lefteqn{ E \Bigl[ \exp{ \Bigl\{  \frac{\m}{\ld_0^{\m}} 
 {\g}_0^2 + \frac{\m}{\ld_4^{\m}} {\g}_4^2 \Bigr\} }  \Bigl|  
    \frac{c_0}{\sqrt{\ld_0^{\m}}}  {\g}_0 -  \frac{c_4}{\sqrt{\ld_4^{\m}}} {\g}_4 = 0 \Bigr]  } \\
& = &  {\sqrt{ 2 \pi (\frac{c_0^2}{\ld_0^{\m}} + \frac{c_4^2}{\ld_4^{\m}})} } 
       \times \int_{-\infty}^{\infty} 
       \exp{ \Bigl\{ - \frac{1}{2} \Bigl[ (1 - 2\m/\ld_0^{\m}) \frac{c_4^2}{\ld_4^{\m}}  
                         + (1 - 2\m/\ld_4) \frac{c_0^2}{\ld_0^{\m}} \Bigr] x^2 \Bigr\}} 
       \frac{dx}{2\pi}    \\
& = & \Bigl(  \frac{\ld_4^{\m} c_0^2 + \ld_0^{\m} c_4^2}{ (\ld_4^{\m} - 2\mu) c_0^2 
        + (\ld_0^{\m} - 2\mu) c_4^2 } \Bigr)^{1/2}. 
\eeqn
Concerning the prefactor in (\ref{equ:expand1}), when restored to a product over the full range,
one over its square 
reads as $P(\ld) \equiv \prod_{\l=0}^{\infty} (1 - \ld /\ld_{\l}^{\m})$ evaluated at $\ld = 2 \m$.  
The behavior of ${\ld}_{\l}^{\m}$ for $\l \up \infty$ shows that the latter is
entire function of order $1/2$.  Also, as $P(\ld)$  vanishes only at the periodic spectrum of $Q_{\m}$, 
one concludes that it is a constant multiple of $\Delta^2(\ld) - 4$.  In other words,
\[
 \Bigl[  \prod_{ 2 \le \l \le \infty, \l \neq 4} \Bigl( 1 - \frac{2 \m}{\ld_{\l}^{\m}} \Bigr)  \Bigr]^{-1/2}
= \sqrt{ \frac{ \Delta^2(0)  - 4 }{ \Delta^2(2 \mu) - 4} }  \times  
  \Bigl[ \Bigl( 1 - \frac{2 \mu }{\ld_0^{\m}} \Bigr)
  \Bigl(  1 - \frac{2 \mu }{\ld_1^{\m}} \Bigr) \Bigl(  1 - \frac{2 \mu }{\ld_4^{\m}} \Bigr) \Bigr]^{1/2}.
\]
The proof is finished.

\vspace{.2cm}

\noindent{\bf Proof of Proposition \ref{lem:g2} }.  The goal is to convert the $CBM$ integral to one over
Brownian bridge paths. First though, the conditionings
$\int_0^1 p = 0$ and $\int_0^1 p \phi_1^{\mu}$ are removed as follows
\beq
\label{equ:second}
\lefteqn{ E_0^0 \Bigl[ \exp{ \Bigl\{ - \frac{1}{2}  \int_0^1 q_{\mu}(x) p^2(x) dx  \Bigr\} } \Bigr]   } \\
& = & P_0^{-1} \Bigl(  \int_0^1 \phi_1^{\m} p   = 0 \Bigr) P_{q_{\m}} \Bigl(  \int_0^1 p = 0, 
   \int_0^1 \phi_1 p = 0 \Bigr) 
      \times CBM \Bigl[  e^{ -\frac{1}{2} \int_0^1 q_{\mu}(x) p^2(x) dx}  \Bigr]. \nonumber  
\eeq
The first factor, while computable, cancels the same object in the 
numerator of (\ref{equ:g2}). 
The second is of a similar nature: under $P_{q_{\m}}$, the variables  $\int_0^1 p(x)dx $ and 
$\int_0^1 \phi_1(x) p(x) dx $ 
are independent Gaussians:
\beq
\label{equ:keep2}         
 P_{q_{\m}} \Bigl(  \int_0^1 p = 0, \int_0^1 \phi_1^{\m} p = 0 \Bigr)  
       & =  &   P \Bigl( \frac{c_0}{\sqrt{\ld_0^{\m}}} \g_0 - \frac{c_4}{\sqrt{\ld_4^{\m}}} \g_4 = 0 \Bigr) 
           P \Bigl( \frac{1}{\sqrt{\ld_1^{\m}}} \g_1 = 0   \Bigr)  \\
  & = &  \frac{1}{\sqrt{ 2 \pi ( {c_0^2}/ \ld_0^{\m} + {c_4^2}/\ld_4^{\m} )}} \times  
 \frac{1}{\sqrt{2 \pi / \ld_1^{\m}}},  \nonumber
\eeq
as advertised.
Turning to the third factor in (\ref{equ:second}), we now unravel the
periodic boundary conditions. From the definition, 
\beq
\label{equ:tiedto}
\lefteqn{ \hspace{-1.7cm}
  CBM \Bigl[  \exp{ \Bigl\{ -\frac{1}{2} \int_0^1  q_{\m}(x) p^2(x) dx  \Bigr\} } \Bigr] 
   =   \int_{-\infty}^{\infty} BM_{00} \Bigl[ \exp{ \Bigl\{ -\frac{1}{2}  
        \int_0^1  q_{\m}(x) ( p(x) + c)^2  dx \Bigr\} } \Bigr]   
          \, \frac{dc}{\sqrt{2 \pi}} }  \\
 & = &  {Z_{q_{\m},0}} \int_{-\infty}^{\infty} \exp{ \Bigl\{ - \frac{1}{2} c^2 
     \int_0^1 q_{\m}(x) dx \Bigr\}  }   
                       E_{q_{\m},0} \Bigl[ \exp{ \Bigl\{ - c \int_0^1 q_{\m} p  
          \Bigr\} } \Bigr] \, \frac{dc}{\sqrt{2 \pi }}. \nonumber
\eeq
Here, in line two, we have renormalized yet again to introduce
the measure $P_{q_{\mu, 0}}$ on tied paths with weight $ Z_{q_{\m},0}
  = BM_{00} [ \exp{ \{ -(1/2) \int_0^1 q_{\mu} p^2 \}} ]$. 
Said differently, $P_{q_{\mu, 0}}$ is the Gaussian measure with inverse covariance operator
${Q}_{\m,0} =  - d^2/ dx^2 + q_{\mu} $ over paths $p(x)$ subject to $p(0) = p(1) = 0$.
As such, the  integrand in (\ref{equ:tiedto}) may be worked out as
\[
\exp{ \Bigl\{ - \frac{1}{2} c^2 
     \int_0^1 q_{\m}(x) dx \Bigr\}  }  
E_{q_{\m}, 0} \Bigl[ \exp{ \Bigl\{ - c \int_0^1 q_{\m}(x) p(x) dx \Bigr\} }
  \Bigr]  
            =  \exp{ \Bigl\{ -\frac{1}{2} c^2  \Bigl[  \int_0^1
  q_{\m}(x) \Bigl( 1 -  ({Q}_{\m,0}^{-1}  q_{\mu})(x) \Bigr) dx  \Bigr]  \Bigr\} }
\]
by simply using the definition of the Green's function ${
 Q}_{\m,0}^{-1}$. Next, with 
 $\psi(x) =  ({\mathfrak G}_0^{-1}  q_{\mu})(x)$ it is immediate that
\[
   \int_0^1 q_{\m}(x) \Bigl( 1 - ({\mathfrak G}_0^{-1} q_{\mu} )(x) \Bigr)  dx 
   = - \int_0^1 \psi^{\pr \pr}(x)  dx 
                          \equiv \int_0^1 \Bigl[ \psi_0^{\pr \pr}(x) + \psi_1^{\pr \pr}(x) \Bigr] dx
\]
in which 
$\psi_0(x)$ and $\psi_1(x)$ are the increasing/decreasing solutions of 
$ \psi^{\prime \prime}(x) = q_{\m}(x) \psi(x)$ 
over $0 < x < 1$ subject to $\psi_0(0) = 0$, $\psi_0(1) = 1$ and $\psi_1(0) = 1$, 
$\psi_1(1) = 0$.  In terms of the 
normalized cosine and sine-like solutions
at $\ld = 0$,  we  have $\psi_0(x) = y_1(x,0) - (y_1(1,0)/y_2(1,0)) y_2(x,0)$,
$\psi_1(x) = y_2(x,0)/y_2(1,0)$, and so
\be
\label{equ:sincos}
   \int_0^1 \Bigl[  \psi_0^{\pr \pr}(x) + \psi_1^{\pr \pr}(x) \Bigr] dx 
   =   
   y_1^{\pr}(1,0) - \frac{ y_1(1,0)}{y_2(1,0)} \Bigl[ y_2^{\pr}(1,0) - 1 \Bigr]
                      + \frac{1}{y_2(1,0)} \Bigl[ y_2^{\pr}(1,0) - 1 \Bigr].            
\ee

Finally, we bring in the  classical computation (see \cite{Sh}),
\[
   Z_{q_{\m} 0} = BM_{00} \Bigl[  \exp{ \Bigl\{ - \frac{1}{2} \int_0^1 q_{\m}(x)  p^2(x) dx \Bigr\} } \Bigr] 
                = \frac{1}{\sqrt{y_2(1,0)}},
\]
which, when combined with (\ref{equ:tiedto}) through (\ref{equ:sincos}), gives
\beq
CBM \Bigl[ \exp{ \Bigl\{ - \frac{1}{2} \int_0^1 q_{\m} (x)  p^2(x) dx  \Bigr\} }  \Bigr] 
& =  & \Bigl[ y_1(1,0) + y_2^{\pr}(1,0)  - 2 \Bigr]^{-1/2}   \\
& =  & \Bigl[  \Bigl( y_1(1/2,0) + y_2^{\pr}(1/2,0) \Bigr)^2  - 4 \Bigr]^{-1/2}  
              \equiv \frac{1}{\sqrt{ \Delta^2(0) - 4}}. \nonumber
\eeq
Here we have made use of the Wronskian identity $1 = y_1 y_2^{\pr} - y_1^{\pr} y_2$ 
in line one and the connection formula
$y_{1,2}(1) = y_{1,2}(1/2) y_{1}(1/2) + y_{1,2}^{\pr}(1/2) y_2(1/2)$ in line two.  
The proof is finished.

\section{Putting it all together: final asymptotics}
\setcounter{num}{7}
\setcounter{equation}{0}
\label{sec:last}

The results through this point are summarized in the following Theorem.

\begin{thm} 
\label{thm:mainish}
For large $\m > 0$,
\be
\label{equ:mainish}
   f(\m) = \sqrt{\frac{2}{\pi}}  A(p_{\m}) R(p_{\m}) Z(\m) \exp{\Bigl[ I_{\m}(p_{\m}) \Bigr]} ( 1 + o(1) )
\ee
with $A(\cdot), R(\cdot), I_{\m}(\cdot)$ and $Z(\mu)$ as defined in 
(\ref{equ:A}), (\ref{equ:defR}), (\ref{equ:Imu}) and (\ref{equ:Z}) respectively.  
Everything on the right hand side is an explicit functional
of $p_{\m} = \sqrt{\m} k \sn(\cdot, k)$  and the simple
spectrum of $Q_{\m}$.
\end{thm}

This is really the main result of the paper. By working out 
of the asymptotics of the individual objects on its right hand side, 
(\ref{equ:mainish}) may be translated to the 
statement given in  Theorem \ref{thm:main}.

To begin, that $I_{\m}(p_{\m}) = - 8/3 {\m}^{3/2} $ up to exponentially small
corrections in $\sqrt{\m}$ has already
been noted.  By similar considerations we find that,
\beq
\label{equ:rr}
  R(p_{\m}) & = & \frac{1}{2}  \m k^2 
 \frac{ \int_0^1 \sn^2(\sqrt{\m} x) \Bigl( \dn^2(\sqrt{\m} x) + k^2 \cn^2( \sqrt{\m} x) \Bigr) dx}
      { \sqrt{ \int_0^1 \cn^2(\sqrt{\m} x) \dn^2(\sqrt{\m} x) dx }}  \\
 & = &  \sqrt{2}  \mu^{3/4}  \frac{\int_{-\infty}^{\infty} \tanh^2(x) \sech^2(x) dx}
                     { \sqrt{    \int_{-\infty}^{\infty} \sech^4(x) dx}} \, 
      \Bigl( 1 + O\Bigl( e^{-{\sqrt{\m}}/4} \Bigr) \Bigr) 
= 
 \sqrt{\frac{2}{3}} {\m}^{3/4} \Bigl( 1 + O\Bigl( e^{-{\sqrt{\m}}/4} \Bigr) \Bigr). \nonumber
\eeq
For $A(p_{\m}) = A_{+}(p_{\m}) A_{-}(p_{\m})$, 
we have first by direct computation: 
\beqn
  A_{+}(p_{\m}) =  \int_0^1 e^{2 k \sqrt{\m} \int_0^x \sn( \sqrt{\m} x^{\pr}) dx^{\pr}} dx
  & = & \frac{1}{\sqrt{(1-k)^2}} 
            \Bigl( \int_0^1 \dn^2(\sqrt{\m} x) + k^2 \cn^2 (\sqrt{\m} x) \Bigr)  dx \\
  & = & \frac{1}{8 \sqrt{\m}} e^{\sqrt{\m}} \times 
           \Bigl( 1 + O\Bigl( e^{-\sqrt{\m}/4} \Bigr) \Bigr). \nonumber
\eeqn
The second half of $A$ responds to
\beqn
A_{-}(p_{\m}) &  =  & \int_0^1 e^{-2 k \sqrt{\m} \int_0^x \sn(\sqrt{\m} x^{\pr}) dx^{\pr}} dx 
  =  \frac{2}{\sqrt{\m}} \int_{0}^{\sqrt{\mu}/4} e^{-2 k \int_0^x \sn(x^{\pr}) dx^{\pr}} dx 
       + O\Bigl(e^{-\sqrt{\m}} \Bigr) \\
 & = &  \frac{2}{\sqrt{\m}} \int_{0}^{\infty} e^{-2 \int_0^x \tanh(x^{\pr}) dx^{\pr}} dx 
            \times \Bigl( 1 + O\Bigl(e^{-\sqrt{\m}/4} \Bigr) \Bigr)
         = \frac{2}{\sqrt{\m}} \times
                  \Bigl(1 +  O\Bigl( e^{-\sqrt{\m}/4} \Bigr) \Bigr) \nonumber,
\eeqn
the integral in the second line being easily computed.

Last we 
turn to the asymptotics of $Z(\m)$.  
Examining the results of Propositions \ref{lem:g1} and \ref{lem:g2}, we see
that we have two objects not  
involving Hill's discriminant:
\[
 \Bigl[  \Bigl( 1 - \frac{2 \m }{\ld_0^{\m}} \Bigr) 
               \Bigl(  1 - \frac{2 \m }{\ld_1^{\m}} \Bigr) 
               \Bigl(  1 - \frac{2 \m }{\ld_4^{\m}} \Bigr) 
\Bigl(  \frac{\ld_4^{\m} c_0^2 + \ld_0^{\m} c_4^2}{ (\ld_4^{\m} - 2\mu) c_0^2 
           + (\ld_0^{\m} - 2\mu) c_4^2 
        } \Bigr)  \Bigr]^{1/2}
 =   \frac{8}{3\sqrt{2}} \m^{1/4} e^{-\sqrt{\m}/2} \times  
         \Bigl( 1 + O \Bigl( \frac{1}{\sqrt{\m}} \Bigr) \Bigr), \nonumber
\]
and 
\be
\label{equ:lastbit}
 \Bigl[  \frac{2 \pi}{\ld_1^{\m}} ( \frac{c_0^2}{\ld_0^{\m}} + \frac{c_4^2}{\ld_4^{\m}} ) \Bigr]^{-1/2}  
          =  \sqrt{\frac{6}{\pi}} \m \times \Bigl( 1 + O \Bigl( \frac{1}{\sqrt{\m}} \Bigr) \Bigr).
\ee
Here we have once again made use $1 - k^2 = 1 - k^2(\m) = 16 e^{-\sqrt{\m}/2} (1 +o (1))$, 
the list (\ref{equ:eigs}),
as well as $c_0(\m) = \sqrt{6} \m^{-1/4} + O(e^{- \sqrt{\m}/4})$ and $c_4(\m) = 1 +  O(e^{- \sqrt{\m}/4})$, 
as pointed out in the proof of Lemma \ref{lem:starmeas}. 
We now make use of Hochstadt's formula to estimate the discriminant for
          large values  of $\m$.

\begin{prop}  
\label{lem:discrim}
It holds
\be
\label{equ:rootdisc}
    \Bigl( \Delta^2( 2 \m) - 4 \Bigr)^{-1/2}  =  e^{ - \sqrt{\mu} } \times  
                             \Bigl(1 + O \Bigl( \frac{1}{\sqrt{\m}} \Bigr) \Bigr)
\ee
as $\m \ra \infty$.
\end{prop}

Gathering one over the right hand side of (\ref{equ:rootdisc})  
together with displays (\ref{equ:rr}) through (\ref{equ:lastbit}) produces 
the form of the result originally stated in Theorem \ref{thm:main}. 

\bigskip

\noindent{\bf Proof of Proposition \ref{lem:discrim} }  Starting from 
$\Delta(\cdot) = 2 \cos \psi( \cdot)$ with $\psi$ defined 
in  (\ref{equ:hoch}) we have that
\beq
\label{equ:lastdis}
 \Delta (2 \mu) 
  & =  &  2 \cos \left( \frac{\sqrt{-1}}{2} \int_{\ld_0^{\m}}^{2\m} 
              \frac{(s - (\ld_1^{\pr})^{\m} ) (s - (\ld_2^{\pr})^{\m} )}{{\sqrt{-(s- \ld_0^{\m}) \cdots
                  (s - \ld_4^{\m})     }   }}  ds \right) \\
  & = &  - 2 \cosh  \left( \frac{\sqrt{\m}}{2}  
             \int_{\ld_1}^{2} (s - \ld_1^{\pr} ) (s - \ld_2^{\pr}) \frac{ ds}{\sqrt{R(s)}} \right) \nonumber.
\eeq
Besides an obvious change of variables, the second line makes use
of  the fact $\psi(\ld_1^{\m})  =  \pi$. 
$R(s)$ is shorthand for  $(s - \ld_0)(s - \ld_1) \times (\ld_2 - s) (\ld_3 - s) (\ld_4 - s)$
which is non-negative for  $s \in [\ld_1,2]$.   
Note that $\ld_0, \dots, \ld_4$, $\ld_1^{\pr}$, and $\ld_2^{\pr}$  now 
refer to the spectral points of the unscaled operator with periodic boundary conditions  
over $[0, 2K = \sqrt{\m}/2]$. While these points still depend on
$\m$, we will show that 
\[
\int_{\ld_1}^2  (s - \ld_1^{\pr} ) (s - \ld_2^{\pr}) \frac{
  ds}{\sqrt{R(s)}}  \ra 2
\]
as $\m  \ra \infty$.  First though we
 must pin down  the behavior of $\ld_1^{\pr} \in [\ld_1, \ld_2]$ and  $\ld_2^{\pr}
\in [\ld_3, \ld_4]$ in that limit. 

Returning to Hochstadt's result we know that 
\be
\label{equ:primes}
   0 = \int_{\ld_1}^{\ld_2} (s - {\ld_1^{\pr}}) ( s - {\ld_2^{\pr}}) \frac{ds}{ \sqrt{R(s)}}
 \ \ \mbox{ and } 
\ \ 
  0 =  \int_{\ld_3}^{\ld_4} (s - {\ld_1^{\pr}}) ( s - {\ld_2^{\pr}}) \frac{ds}{\sqrt{R(s)}}
\ee
from which it may be immediately inferred that $|\ld_1^{\pr} -  \ld_1|$ and $|\ld_2^{\pr} - \ld_3|$ 
tend to zero 
as $\m \ra \infty$.
(Recall that, as $\m \ra \infty$, $\ld_0 \simeq \ld_1 = 2 + O(e^{-\sqrt{\m}/2})$ and 
$\ld_2 \simeq \ld_3 = 5 +O(e^{-\sqrt{\m}/2})$.) 
Indeed, if $\ld_2^{\pr}$ were to remain greater than $ \ld_3 + \delta $ for a fixed 
$\delta > 0$, the second equality
in (\ref{equ:primes}) would require that $\ld_1^{\pr} \ra 5$ as $\mu \ra \infty$.  However, with
$\ld_1^{\pr} \simeq 5$ there would be nothing to balance the singularity
at the lower limit of the integral in the first equality. 
Given now that $\ld_2^{\pr} \ra 5$, it must be that $\ld_1^{\pr} \ra 2$ 
in order that the first integral remains finite.  

For sharper information 
it is convenient to introduce the (small) parameter $\ep \equiv 1 - k^2$, whereupon
\beq
\label{equ:epses}
  \ld_0 = 2 - \ep - 3/4 \,\ep^2 - O({\ep^3}), & &   \ld_1 = 2 - \ep,  \\
\ld_2 = 5 -  4 \ep , \ \ \ \  \ \ld_3 = 5 - \ep,  &  &  \ld_4 = 6 - 4 \ep -  O(\ep^2). \nonumber
\eeq
We further denote $\ld_1^{\pr} = \ld_1 + \d_1(\ep)$ and $\ld_2^{\pr} = \ld_3 + \d_2(\ep)$.
The asymptotics of $\d_2(\ep)$ as $\ep \ra 0$ may now be determined from the second relation 
in (\ref{equ:primes}).  Note first that since  $(s- \ld_1^{\pr})/\sqrt{ (s - \ld_0) (s- \ld_1) } = 1 + o(1)$ 
for $s \in [\ld_3, \ld_4]$ and $\ep \ra  0$, it suffices to consider
\[ 
    0 = \int_{\ld_3}^{\ld_4}  \frac{ ( s - \ld_3 - \d_2(\ep) ) }{ \sqrt{ ( s - \ld_2) ( s - \ld_3) ( \ld_4 - s)}} ds
\]
in order to gather the leading order behavior of $\d_2(\ep)$.  Next, shifting variables and employing
(\ref{equ:epses}), this last equality is the same as
\[
    \int_0^1 \sqrt{ \frac{s}{s+ 3 \ep}} \frac{ds}{\sqrt{1-s}} = 
   \d_2(\ep)  \int_0^1 \frac{ ds}{\sqrt{ s ( s + 3 \ep) (1 -s) }}
\] 
us to negligible errors.
From here one quickly concludes that $\d_2(\ep) = 2 ( \log |\ep| )^{-1} (1 + o(1))$. Armed 
with this information, similar considerations brought to bear in the first relation in
(\ref{equ:primes}) lead to the estimate $\d_1(\ep) = 4 ( \log |\ep| )^{-1} (1 + o(1))$.

Now we may return to the main integral at hand:
\begin{eqnarray*}
\lefteqn{
\int_{\ld_1}^{2} (s - \ld_1^{\pr} ) (s - \ld_2^{\pr}) \frac{ ds}{\sqrt{R(s)}}
   =  \int_0^{\ep}  ( s - \d_1(\ep) ) (s - (\ld_2^{\pr} - \ld_1) ) \frac{ ds}{\sqrt{R(s + \ld_1)}}}  \\
& = & \Bigl(  \int_0^{1}  \frac{ (\d_1(\ep) - \ep s )}
                                 { \sqrt{s ( s + 3/4 \ep + o(\ep))}}  
                   \frac{ (3 + \d_2(\ep) - \ep s ) ds}{ \sqrt{(3 - \ep s) ( 3 - \ep  - \ep s) (4 - s) }}  
                  \Bigr)   ( 1 + o(1)) \\
& = &   \Bigl(   \d_1( \ep) 
              \int_0^{1}  \frac{ds}
                                 { \sqrt{s ( s + 3/4 \, \ep)}}  - \ep \int_0^1  
                    \sqrt{ \frac{ s }{s + 3/4 \, \ep}} ds  
                                 \Bigr) \, \Bigl( \frac{1}{2} + O(\d_2(\ep)) \Bigr)  \\ 
& = & \frac{1}{2} \Bigl( \d_1(\ep) \log ( {1}/{\ep})  + O ( \d_1(\ep) ) \Bigr) \Bigl( 1 + O(\d_2(\ep)) \Bigr)\\
 & = & 2 +  O \Bigl( \log^{-1}( 1/\ep)  \Bigr).
\end{eqnarray*}
That is, the right hand side of (\ref{equ:lastdis}) equals $(-2)$ times  $ \cosh( \sqrt{\m} + O(1))$.
The proof is finished.

\section{Appendix}
\setcounter{num}{8}
\setcounter{equation}{0}

\subsection{On the rate function $I$}

In this first section of the appendix we revisit the minimization problem (\ref{equ:Imu}),
proving a few technicalities used in the proof of Theorem \ref{thm:rate}.

\begin{prop}  Consider the Euler-Lagrange equation corresponding to the minimizer
of $I(f; a) $:
\be
\label{equ:EL2}
     \frac{1}{2} ( f^{\prime})^2 = \frac{1}{2} f^4 - f^2 - \a f + \frac{1}{2} \beta,
\ee
recall (\ref{equ:EL1}).  Here $\a$ is the multiplier and $\b$ a constant of integration. As
$a \ra \infty$, $\a \ra 0$ and $\b \ra 1$.
\end{prop}

It is natural to suppose that $\a \equiv 0$ for all $ a $
sufficiently large, but we were unable to prove this. 

\bigskip

\noindent{\bf Proof } We first show that $|\b - 1| = O(1/a)$.   Integrating (\ref{equ:EL2}) 
and recalling that we are working in the space of mean-zero functions we find
\be
\label{equ:ap1}
   \int_{-a}^{a} ( d f^*/dx)^2 = \int_{-a}^a [ 1 - (f^*)^2]^2 + 2 (\b -1) a
\ee
for any minimizer $f^* = f_a^*$. It follows that
\[
  \frac{8}{3} \ge I^*(a) = \int_{-a}^a [ 1 - (f^*)^2]^2 +  (\b - 1) a \ge ( \b - 1) a, 
\]
which gives us one direction.  For the opposite inequality, solve instead for
$  \int [ 1 - (f^*)^2]^2 $ in (\ref{equ:ap1}) and substitute that into $I(f^*;a) = I^*(a)$
to yield
\[
  \frac{8}{3} \ge I^*(a) = \int_{-a}^a ( d f^*/dx)^2  +  (1 - \b) a \ge ( 1 - \b) a. 
\]

Turning to the convergence of $\a$ to zero, integrating the preliminary Euler equation,
$f^{\pr \pr} = 2 f^3 - 2 f - \a$,  implies
$
   \a  = - 2 \int_{-a}^a (f^*)^3,
$
and so, since $f^*$ and $-f^*$ both minimize, it may from here on be assumed that $\a \le 0$.
The important point is to see that $|| f^*||_{\infty} \le 1$ and that $|| f^*||_{\infty} 
\ra 1$ as $a \ra \infty$. First let $\theta$ be a maximum of $f^*$.  At any point it is attained
(\ref{equ:EL2}) reads as $0 = (1-\t)^2 + 2 |\a| \t + (\b - 1)$ which shows that
\[
     2 \t |\a| \le 1 - \b.
\]
Now assume that $\t > 1$.  Then there is a point at which $f^* = 1$ and
\[
   0 \le 2 |\a| + \b - 1 < 2 |\a| \t + \b  - 1 \le 0 
\]
by (\ref{equ:EL2}) and the previous display. But this is a contradiction.  A similar
argument explains why $f^*$ is everywhere greater than $-1$.  Finally denoting 
$\eta  = f^*(x^*)$ where $| f^*(x^*)| = ||f^*||_{\infty}$ it is plain that $\eta \ra 1$
if $I^*(a)$ is to remain bounded.  And, again from (\ref{equ:EL2}) now at $x^*$, we
have that
\[
   0 = ( 1 - \eta^2)^2 + 2 |\a| \eta + (\b - 1) 
\] 
from which it follows that $\a \ra 0.$

\subsection{Two changes of measure}

Here we provide the proof of the result which allowed us to remove the degeneracy, 
Lemma \ref{l:StatZeros},  as well as that of the Cameron-Martin formula for $P_0^0$ invoked
in the proof of Proposition \ref{p:shift}.   

As the reader may have observed, Lemma \ref{l:StatZeros} is really a type
of Rice formula (see \cite{Rice}, \cite{Kac}, or \cite{Lead} for a more recent account).   
A proof is provided as we were unable to locate the form of the result needed here.  Our
proof relies on finite dimensional approximation, and in that direction we
first prepare the following.

\begin{lemma}
\label{l:finitecond}
Suppose that $X(\cdot)$ is a stationary process, periodic of period one,
defined on the lattice  
$\{  k /2^n \} $ and taking values in $\Z/n$.
Assume also that $X$ has at least one zero with probability $1$. If ${ F}$ is a 
functional that is invariant under translations, then
\[
E[{ F}(X)] = 2^n E \Bigl[ { F}(X) N^{-1} \Bigl| X(0)=0 \Bigr] P \Bigl( X(0)
= 0 \Bigr)
\]
where $N=N(X) = $ the number of zeros of $X$.
\end{lemma}

\noindent{\bf Proof } Notice the simple decomposition
\begin{eqnarray*}
E[{ F}(X),N=1] & = & \sum_{k=1}^{2^n} E \Bigl[ { F}(X), N=1, T_1=k/2^n \Bigr] \\
               & = & 2^n E \Bigl[ F (X), N=1, X(T_1) = 0 \Bigr].
\end{eqnarray*}
Here $T_1$ is the (only) zero of $X$, and we have used the rotation
invariance of $F(X)$ in line two.  A more general version of this is
\begin{eqnarray*}
m \, E \Bigl[ { F}(X),N=m \Bigr]  & = & \sum_{i=1}^m \sum_{k=1}^{2^n} 
E \Bigr[ { F}(X), N=m, T_i=k/2^n \Bigr]\\
& = & \sum_{k=1}^{2^n} E \Bigl[ { F}(X), N=m,\bigcup_{i=1}^m  \{ T_i=k/2^n
  \} \Bigr]\\
& = & \sum_{k=1}^{2^n} E \Bigl[ { F}(X), N=m,X({k/2^n})=0 \Bigr]
 =  2^n  E \Bigl[ { F}(X), N=m, X(0) = 0 \Bigr] .
\end{eqnarray*}
Now divide 
both sides of  this equality by $m$ and sum to conclude the proof of the lemma.

\bigskip

\noindent{\bf Proof of Lemma \ref{l:StatZeros} }
If now $X$ is smooth stationary periodic process over the line we 
define the discrete process $Y^n$ as  taking the values $n^{-1} [n
  X({k/2^n})]$ at $1/2^n, 2/2^n, \dots$.  
Applying the previous result to $Y^n$ produces
\[
E \Bigl[{F}(Y^n)\Bigr] =   
2^n E \Bigr[{ F}(Y^n)  N(Y^n)^{-1} \Bigl| 0\leq X(0) < 1/n \Bigr] P \Bigl( 0
\leq X(0) < 1/n \Bigr)
\]
Now notice that the number of zeros 
of $Y^n$ around a zero of $X$ is approximately $2^n/n\cdot|X^{\prime}(z)|^{-1}$. In fact
\[
 N(Y^n) \simeq \frac{2^n}{n} \sum_{z\in Z} |X^{\prime}(z)|^{-1}
\]
can be used in the above formula to produce
\[
E \Bigl[ { F} (Y^n) \Bigr] \simeq 
E \Bigl[ {F}(Y^n) \Bigl(  \sum_{z\in Z} |X^{\prime}(z)|^{-1}  \Bigr)^{-1} 
\Bigl|  0\leq X(0) <1/n \Bigr] \left( 
\frac{P(0\leq X(0) <1/n)}{1/n} \right) .
\]
This expression has the required limit when $n\uparrow \infty$.  
The proof is finished.

\bigskip

Finally we have:

\bigskip

\begin{lemma}
\label{l:CMart}
Let $\tilde{E}$ denote the $CBM$ conditioned on both
$\int_0^1 p = 0 $ and $\int_0^1 \phi p = 0$ for some continuous,
mean-zero $\phi$.  Let also
 $\varphi$ be twice continuously differentiable and satisfy 
$\int_0^1 \varphi = 0 $ as well as  
$ \int_0^1 \phi  \varphi = 0$. 
Then
\[
    {\tilde E} \left[ F{\left( p\right)} \right] =
     {\tilde E}  \left[ F{\left( p+\varphi \right)}
\exp \left\{ \int_{0}^{1}\varphi ^{\prime \prime }p-\frac{1}{2}
\int_{0}^{1}\left| \varphi ^{\prime }\right|^{2} \right\} \right].
\]
for all bounded measurable functions $F$ of the path.
\end{lemma} 

\noindent{\bf Proof }
The $\tilde{E}$ integral is expressed as an integral with respect to
the Brownian bridge measure:   
\beqn
\tilde{E} \Bigl[ F(p) \Bigr] 
 & = & BM_{00} \left[ F \left( p - \int_0^1 p \right) \ \Bigl| \ \int_0^1 p \phi = 0 \right]  \\
 & = & \lim_{h \downarrow 0} \lim_{\ep \downarrow 0} \frac{1}{B_{h,\ep}} 
      BM_0 \left[  F \left( p - \int_0^1 p \right),  \int_0^1 p \phi \in [0, h], 
                  \ p(1) \in [0, \ep] \right]
\eeqn
where $B_{h,\ep}$ 
$ = BM_0 [   \int_0^1 p \phi \in [0, h], p(1) \in [0, \ep] ]$.
The usual Cameron-Martin formula may now be applied to the free Brownian integral
with the result that
\beq
\label{equ:CM1}
\lefteqn{BM_0 \left[  F \left( p - \int_0^1 p \right),  
              \int_0^1 p \phi \in [0, h], p(1) \in [0, \ep] \right]} \\
& = &  BM_c \left[ F \left( p - \int_0^1 p + \varphi \right) 
                   \exp\Bigl\{ - \int_0^1 \varphi^{\pr} d p - \frac{1}{2} \int_0^1 |\varphi^{\pr}|^2 \Bigr\}
                     \int_0^1 p \phi \in [0, h], \ p(1) \in [c, c + \ep] \right],
                    \nonumber
\eeq
in which $ c = - \varphi(0)$ and use has been made of $\int_0^1 \varphi = 0$ and 
$\int_0^1 \varphi \phi = 0$.  Next It${\hat{\mbox{o}}}$'s formula 
and the fact that $ \int_0^1 \varphi^{\pr \pr} = 0 $,  allows us to write
\[
- \int_0^1 \varphi^{\pr}(x) d p(x)  = \int_0^1 \varphi^{\pr \pr}(x) 
         \Bigl( p(x) - \int_0^1 p(x^{\pr}) dx^{\pr} \Bigr) dx   
                                  + \varphi^{\pr}(0) \Bigl( p(1) - p(0) \Bigr).  
\]
Invoking this move as well as the shift from $p(x)$ under 
$BM_c$ to $p(x) + c$ under $BM_0$, the right hand side of (\ref{equ:CM1}) 
then reads
\beqn
BM_0 \left[ F \left( p - \int_0^1 p + \varphi \right) 
 \exp \Bigl\{ - \int_0^1 \varphi^{\pr \pr } \Bigl( p - \int_0^1 p \Bigr) 
                        - \frac{1}{2} \int_0^1 |\varphi^{\pr}|^2 \Bigr\} \right. \\
 \left. \hspace{1cm}  \exp \Bigl\{ - \phi^{\pr}(0) p(1) \Bigr\},         
                \int_0^1 p \phi \in [0, h], \ p(1) \in [0, 0 + \ep] \right].                
\eeqn
Dividing this object by $B_{h, \ep}$ and performing the limits $\ep \downarrow 0$
and $h \downarrow 0$ completes the proof.

\bigskip

\noindent{\bf Acknowledgments}  Deepest thanks to H.P. McKean for having
introduced S.C. and B.R. to this circle of ideas. B.R.
would also like to acknowledge the support of the NSF grant DMS-9883320.


\begin{thebibliography}{10}


\bibitem{A}
{\sc Adler, R.}
{\it An Introduction to Continuity, Extrema, and Related Topics
for General Gaussian Processes},
IMS Lecture Notes-Monograph Series, No. 12. California, 1990.



\bibitem{GB}
{\sc Ben Arous, G.; Deuschel, J-D; Stroock, D.W. }
Precise asymptotics in 
large deviations.
{\it Bull. Sci. Math} {\bf 117} (1993), no. 1, 107-124.


\bibitem{Bol}
{\sc Bolthausen, E.}
Laplace approximations for sums of independent random vectors. II.
Degenerate maxima and manifolds of maxima.
{\it Probab. Theory Related Fields} {\bf 76} (1987), 167-206.


\bibitem{BF}
{\sc Byrd, P.; Friedman, M.D.}
Handbook of Elliptic Integrals,
Springer-Verlag, Berlin-G\"ottingen-Heidelberg, 1954.

\bibitem{CM}
{\sc Cambronero, S;  McKean, H.P.}
The Ground State Eigenvalue of Hill's Equation
with White Noise Potential.
{\it Comm. Pure Appl.
Math.} {\bf 52} (1999), 1277-1294. 


\bibitem{ER}
{\sc Ellis, R.S.; Rosen, J. S.}
Asymptotic Analysis of Gaussian Integrals, II:
Manifolds of Minimum Points.
{\it Comm. Math Phys.}
{\bf 82} (1981), 155-181. 

\bibitem{FL}
{\sc Frisch, H.L.; Lloyd, S.P.}
Electron levels in a one-dimensional lattice.
{\it Phys. Rev.} {\bf 120} (1960), no. 4, 1175-1189.


\bibitem{FM}
{\sc Fukushima, M.; Nakao, S.}
On the spectra of the Schr\"odinger operator with
a white Gaussian noise potential.
{\it Z. Wahr. und Verw. Gabiete}
{\bf 37} (1976/77), no. 3, 267-274

\bibitem{Hlp}
{\sc  Halperin, B.I.}
Green's functions for a particle in a one-dimensional random
potential.
{\it Phys. Rev. (2)} {\bf 139} (1965), A104-A117.

\bibitem{Hoch}
{\sc Hochstadt, H.}
Functiontheoretic properties of the discriminant of
Hill's equation.
{\it Math. Zeitchr.} {\bf 82} (1963), 237-242.

\bibitem{Hoch1}
{\sc Hochstadt, H.}
Asymptotic estimates for the Sturm-Liouville spectrum.
{\it Comm. Pure. Appl. Math} {\bf 14} (1961), 749-764.


\bibitem{I} 
{\sc Ince, E.L.} 
The periodic Lam{\'e} functions.
{\it Proc. Roy. Soc. Edinburgh.} {\bf 60} (1940), 47-63.


\bibitem{Kac}
{\sc Kac, M.}
On the average number of real roots of a random algebraic equation.
{\it Bull. Amer. Math. Soc.} {\bf 49} (1943), 314-320.


\bibitem{KS1}
{\sc Kusuoko, S.; Stroock, D. W.}
Precise asymptotics of certain Wiener functionals.
{\it J. Funct. Anal.}
{\bf 99} (1991),  1-74.


\bibitem{KS2}
{\sc Kusuoko, S.; Stroock, D. W.}
Asymptotics of certain Wiener functionals with degenerate extrema.
{\it Comm. Pure. Appl. Math.}
{\bf 47} (1994), 477-501.


\bibitem{Lead}
{\sc Leadbetter, M. R.; Spaniolo, G.V.}
Reflections on Rice's formulae for level crossing $-$ history, extensions
and use.
{\it Aust. N. Z. J. Stat.} {\bf 46} no. 1  (2004) , 173-180.


\bibitem{LGP}
{\sc Lifshits, I.M.; Gredeskul, S. A.; Pastur, L.A.}
{\it Introduction to the theory of disordered systems}.
J. Wiley \& Sons, New York, 1988.

\bibitem{MW}
{\sc Magnus, W; Winkler, S.} 
{\it Hill's Equation.}
Interscience Tracts in Pure and Applied Mathematics,
No. 20. Interscience, J. Wiley \& Sons. 
New York-London-Sydney, 1996.

\bibitem{Mrk}
{\sc Merkl, F.}
Quenched asymptotics of the ground state energy
of random Schr\"odinger operators
with scaled Gibbsian potentials.
{\it Probab. Theory Relat. Fields}
{\bf 126} (2003) 307-338.

\bibitem{M1}
{\sc McKean, H.P.}
A limit law for the groundstate of Hill's equation.
{\it J. Stat. Phys.} {\bf 74} (1994), nos. 5-6 1227-1232.

\bibitem{MvM} 
{\sc McKean, H.P.; van Moerbeke, P.} 
The Spectrum of Hill's Equation.
{\it Inventiones Math.} {\bf 30} (1975), 217-274.


\bibitem{Rice} 
{\sc Rice, S. O.}
Mathematical analysis of random noise.
{\it Bell Sys. Tech. J.} {\bf 23} (1944), 292-232.


\bibitem{Sh}
{\sc Shepp, L.A.}
Radon-Nikodym Derivatives of Gaussian Measures.
{\it Ann. Math. Statist.} {\bf 37} (1966), 321-354.


\bibitem{Sn}
{\sc Snitzman, A.-S.}
Brownian motion, obstacles and random media.
Springer Monographs in Mathematics, Berlin-Heidelberg, 1998.


\bibitem{Var}
{\sc Varadhan, S.R.S.}
Lectures on diffusion problems and partial differential equations,
Tata Institute Lectures on Mathematics and Physics, {\bf 64}, 1980.

\end{thebibliography}
\end{document}